\documentclass{amsart}

\usepackage{amsthm}
\usepackage{amssymb}
\usepackage{amsmath}
\usepackage{amscd}

\newtheorem {Theorem}   {Theorem}

\numberwithin{Theorem}{section}

\newtheorem {Lemma}[Theorem]    {Lemma}         
\newtheorem {Proposition}[Theorem]{Proposition}  
\newtheorem {Conjecture}[Theorem]{Conjecture}  
\theoremstyle{definition}
\newtheorem{Definition}[Theorem]{Definition}
\theoremstyle{remark}
\newtheorem{Remark}[Theorem]{Remark}
\newtheorem{Example}[Theorem]{Example}

\newtheorem {Corollary}[Theorem]{Corollary}

\def    \R      {{\mathbb R}}

\def    \C      {{\mathbb C}}

\def    \Z      {{\mathbb Z}}

\def    \LL     {{\mathbb L}}
\def    \L     {{\mathbb L}}

\def    \A     {{\mathbb A}}
\def    \G     {{\mathbb G}}

\def    \su     {{\mathfrak su}}     
\def    \u     {{\mathfrak u}}     
\def    \gl     {{\mathfrak gl}}

\def    \calF   {{\mathcal F}}

\def    \calV     {{\mathcal V}}

\def    \df     {\Omega^1}
\def    \sec     {\Gamma(E)}

\def    \tr    {{\operatorname{tr}}}
\def    \var    {{\operatorname{var}}}

\def    \eps   {\epsilon}

\def    \ol     {\overline}

\def    \g      {{\mathfrak g}}

\def    \h      {{\mathfrak h}}

\def    \ssminus        {{\smallsetminus}}

\def    \pt     {\mathit{pt}}
\def    \p      {\partial}
\def    \Vect   {{\operatorname{Vect}}}

\def    \id     {\mathit{id}}

\def    \ssminus        {\smallsetminus}

\def    \SU     {\operatorname{SU}}
\def    \U      {\operatorname{U}}

\def    \GL     {\operatorname{GL}}

\def    \Vect   {\operatorname{Vect}}

\def	\ol	{\overline}

\def    \Rep   {\mathit{Rep}}
\def    \hom    {{\operatorname{hom}}}
\def    \defo    {{\operatorname{def}}}

\def    \modd   {\operatorname{mod}}

\def    \Lie    {\operatorname{Lie}}

\def    \ol     {\overline}

\newcommand{\labell}[1] {\label{#1}}
\newcommand    {\comment}[1] {}

\begin{document}

%%%%%%%%%%%%%%%%%%%%%%%%%%%%%%
%   TEXT FORMATTING

\setlength{\smallskipamount}{6pt}
\setlength{\medskipamount}{10pt}
\setlength{\bigskipamount}{16pt}

%%%%%%%%%%%%%%%%%%%%%%%%%%

%%%%%%%%%%%%%%%%%%%%%%%%%%

%%%%%%%%%%%           BEGINNING OF  TEXT

%%%%%%%%%%%%%%%%%%%%%%%%%%

\title[Grothendieck Groups of Poisson Vector Bundles]
{Grothendieck Groups of Poisson Vector Bundles}

\author[Viktor L. Ginzburg]{Viktor L. Ginzburg}

\address{Department of Mathematics, UC Santa Cruz, Santa Cruz, 
CA 95064, USA}
\email{ginzburg@math.ucsc.edu}

\date{\today}

\thanks{The work was partially supported by the NSF and by the 
faculty research funds of the UC Santa Cruz.}

\begin{abstract}
A new invariant of Poisson manifolds, a Poisson $K$-ring, is introduced.
Hypothetically, this invariant is more tractable than such invariants
as Poisson (co)homology. A version of this invariant is also defined
for arbitrary algebroids. Basic properties of the Poisson $K$-ring are
proved and the Poisson $K$-rings are calculated for a number of examples.
In particular, for the zero Poisson structure the
$K$-ring is the ordinary $K^0$-ring of the manifold and for
the dual space to a Lie algebra the $K$-ring is the ring of virtual 
representations of the Lie algebra. It is also shown that the $K$-ring
is an invariant of Morita equivalence. Moreover, the $K$-ring is a functor
on a category, the weak Morita category, which generalizes
the notion of Morita equivalence of Poisson manifolds.
\end{abstract}

\maketitle

\section{Introduction}
\labell{sec:intro}

Poisson manifolds have a wide variety of invariants, almost all
of which, with few notable exceptions, are extremely hard to compute.
For instance, Poisson (co)homology and its various ``derivatives''
are not calculated explicitly even for dual spaces to semisimple
Lie algebras. (In fact, it is not even clear in what terms the result
should be expressed.) The goal of this paper is to introduce a more
tractable invariant of Poisson manifolds, a Poisson $K$-ring.
Alternatively, this paper is devoted to the study 
of representations of the algebroid of differential 
forms on a Poisson manifold, albeit the study is oriented toward
a very specific and narrow objective. This work fits in the subject
of Poisson topology which we loosely define by simply listing some of the
relevant publications: 
\cite{Be,crainic,Fe1,ELW,GG,gi-lu:morita,we:local,we:modular,xu}.

\subsection{Poisson $K$-theory}
Poisson geometry occupies an intermediate place between geometry 
of smooth manifolds and non-commutative geometry. The $K$-theory,
the most robust and algebraic of all geometric ``cohomology theories'',
is fairly easy to define in either smooth and non-commutative settings,
and it should also have an analogue in the context of Poisson geometry.

Leaving for a moment the question of its construction aside, we note that
it is reasonable to expect a ``correct'' Poisson $K^0$-functor $K_\pi$ 
to satisfy the following three conditions:
\begin{enumerate}

\item[(K1)] When $P$ is a Poisson manifold with zero Poisson structure, the
Poisson $K$-ring of $P$ is equal to the ordinary topological 
$K$-ring of $P$,
i.e., $K_\pi(P)=K(P)$. Moreover, $K_\pi(P)=K(B)$, if $P=B\times F$
with $\pi=0\oplus\omega^{-1}$, where 
$B$ is a manifold with zero Poisson structure and $(F,\omega)$ is a
simply connected symplectic manifold.\footnote{Here we adopt the standard
topological convention that a simply connected space is 
by definition connected. We also consider only the $K^0$-functor: 
$K(B)=K^0(B)$ and $K(P)=K^0(P)$.}  

\item[(K2)] When $P$ is symplectic with a sufficiently nice (e.g., finite)
fundamental
group, $K_\pi(P)=R(\pi_1(P))$, where $R$ denotes the ring of representations.
In general, in the symplectic case, $K_\pi(P)$ should lie in the range between
$R(\pi_1(P))$ and the Grothendieck ring generated by the semi-ring of 
path-connected 
components of $\Rep(\pi_1(P))$. (Here $\Rep$ stands for the space of 
representations.)

\item[(K3)] Let $\g$ be a Lie algebra.
The ring $K_\pi(\g^*)$ should reflect representation properties of 
$\g$ and one should have $K_\pi(\g^*)=R(\g)$, when $\g$ is sufficiently nice,
e.g., semisimple. 
In general, $K_\pi(\g^*)$ should lie in the range between
$R(\g)$ and the Grothendieck ring generated by the path-connected components
of $\Rep(\g)$. 
\end{enumerate}

These three requirements deserve some discussion. 

Note that the condition (K1), which we accept in this paper
unconditionally, is incompatible with the requirement that
$K_\pi(\g^*)=R(\g)$ for all $\g$. Indeed, for $\g=\R$, we have
$K_\pi(\R)=\Z$ by (K1), whereas $R(\R)$ is a free group
generated by $\R$. This simple observation is indicative of the difficulties 
arising in the definition of a Poisson $K$-ring and is the main reason 
that the third condition is stated somewhat vaguely.

Furthermore, we emphasize that in (K2) we do \emph{not}
require that $K_\pi(P)=K(B)$, whenever the symplectic foliation
of $P$ is a fibration over $B$ with simply connected fibers. In fact, Poisson
$K$-rings defined in this paper are sufficiently sensitive to detect
whether or not the Poisson structure on $P$ varies from fiber to fiber.

The similarity between the requirements (K2) and (K3) arises from
one of the main principles of Poisson geometry. In this connection we simply 
quote Alan Weinstein, \cite{we:poisson}: \emph{It is tempting to think of the 
symplectic manifold $N$ as the ``dual of the Lie algebra of $\pi_1(N)$''.} 
Conversely, a simply connected group with Lie algebra $\g$ can be viewed
as the ``Poisson fundamental group of $\g^*$''. In fact, $K$-rings
considered in this paper are defined to be ``the rings of virtual 
representations of the Lie algebra whose dual is $P$''.

A $K$-theory need not, of course, arise as the Grothendieck group
of a certain class of vector bundles. This is the case, with
example, for the $K$-theory of foliated manifolds defined in the context
of non-commutative geometry, \cite{connes-survey,connes-book}. This suggests
that a ``correct'' Poisson $K$-ring should also be constructed using 
non-commutative geometry techniques. However, this idea encounters a number 
of technical problems.

Here we take a more naive approach, introduce 
Poisson vector bundles,  and define two versions of a Poisson 
$K$-ring as the Grothendieck rings for suitable classes of vector
bundles and their equivalences. As a consequence, the resulting $K$-rings
cannot be expected to give meaningful results for Poisson manifolds with
severely non-Hausdorff spaces of symplectic leaves. 
On the other hand, for manifolds
whose spaces of leaves are only mildly non-Hausdorff, these rings
appear to be ``correct''. (This is the case, for example, for $\g^*$, when
$\g$ is semisimple.) In brief, the
$K$-rings introduced in this paper are related to the ``correct'' one,
which may fail to exist, in the same way as the classical $K$-theory is
related to the non-commutative one. 

The objects we refer to as Poisson vector bundles are well known
in a different context as representations of algebroids. Already in the
simplest examples, they form too large a semi-ring to produce a
Grothendieck ring satisfying (K1)-(K3). We deal with this problem in two
ways: by imposing an extra condition on representations and by relaxing
the equivalence relation. The resulting semi-rings are then suitable for
our purposes. 

\subsection{Organization of the paper; conventions} 
In Section \ref{sec:pvb} we define
Poisson vector bundles and $K$-rings, prove their basic properties, 
and examine a number of examples. No reference to algebroids and groupoids
and their representations is made here. In Section \ref{sec:algebroids},
we put the constructions of Section \ref{sec:pvb} in the context of
representations of algebroids and groupoids. We also recall the
definition of secondary characteristic classes due to Crainic, 
\cite{crainic}, and establish their properties important for what follows.
Section \ref{sec:g} is entirely devoted to the calculation 
of the Poisson $K$-rings for the dual space to a Lie algebra. In
Section \ref{sec:mvb} we prove Morita invariance of Poisson $K$-rings.
In Sections \ref{sec:m-w} and \ref{sec:proof} we partially generalize 
this result to representations of algebroids and 
introduce a new category (the weak Morita category) on which algebroid 
cohomology and one of our $K$-rings are functors.

A few meaningful questions and calculations are entirely left out in the
present account. Among these are the calculations of $K$-rings for 
Bruhat--Poisson structures and for irrational affine Poisson structures
on tori. In fact,
getting explicit answers in these examples would be the ultimate validity
test for $K$-rings as tractable invariants of Poisson manifolds. Neither
do we touch the problem of finding an analogue of Bott periodicity or
defining $K^n$-functors. We only consider Poisson and algebroid versions
of $K^0$. 
Note also that our $K$-rings are defined using a very narrow
class of vector bundles over a Poisson manifolds or representations
of an algebroid. There are some indications that more general objects, 
representations up to homotopy, are at least as natural in the context of
Poisson or algebroid geometry as genuine representations; 
see \cite{Fe1,Fe2,ELW}.
We do not work with this bigger class, for the resulting Grothendieck rings
appear to be more difficult to calculate even in simple examples.

Much of what follows holds simultaneously for real or orthogonal or 
complex or unitary vector bundles with obvious modifications. Hence, we only
specify the class of vector bundles when the distinction is relevant. 
The reader may view the case of complex vector bundles as the 
``default setting'' throughout this paper.

\subsection*{Acknowledgments.} The author is grateful to Rui Loja Fernandes 
and Alan Weinstein for useful discussions. He would also like to express
his gratitude to the EPFL for its hospitality during the period when some
parts of this paper were written.

\section{Poisson Vector Bundles}
\labell{sec:pvb}
\subsection{Definitions}
Let $P$ be a Poisson manifold with Poisson structure $\pi$.

\begin{Definition}
\labell{def:main}
A \emph{Poisson vector bundle} $E\to P$ over $P$ is a vector bundle 
together with a bracket 
$\{~,~\}\colon C^\infty(P)\times \Gamma(E)\to\Gamma(E)$ which turns
$\Gamma(E)$ into a Lie algebra module over $C^\infty(P)$ and such
that the following two Leibniz identities hold
\begin{enumerate}

\item[(L1)] $\{f, hs\}=\{f,h\}s+h\{f,s\}$ and

\item[(L2)] $\{fh, s\}=f\{h,s\}s+h\{f,s\}$ 
\end{enumerate}
for all $f$ and $h$ in $C^\infty(P)$ and all $s\in \Gamma(E)$.
\end{Definition}

Isomorphisms of Poisson vector bundles are defined as isomorphisms of
ordinary vector bundles preserving the Poisson vector bundle structure.
We denote the semi-ring of isomorphism classes of
Poisson vector bundles over $P$ by $\Vect_\pi(P)$. 
The pull-back of a Poisson vector bundle under a Poisson map does not
naturally receive the structure of a Poisson vector bundle, i.e., 
$\Vect_\pi$ is not a functor for Poisson maps. In particular,
the restriction of a Poisson vector bundle to a symplectic leaf does not 
have, in general, a natural structure of a Poisson vector bundle.

The two Leibniz identities in the definition of Poisson vector bundles
play very different roles. The identity (L1) appears to be absolutely
necessary and simply expresses the fact that $\Gamma(E)$ is a Poisson
module over $C^\infty(P)$. The second identity (L2) is more technical as 
the following obvious observation indicates:

\begin{Proposition}
\labell{prop:Leibn2}
The second Leibniz identity (L2) is equivalent
to requiring that $\{f,s\}(p)$, for all points $p\in P$, depends only on 
the differential $df(p)$ (and, by (L1), on the 1-jet of $s$ at $p$).
The second Leibniz identity is also equivalent to the following 
two conditions:
\begin{equation}
\labell{eq:Leibn21}
\{1, s\}=0\quad\text{for all}\quad s\in \Gamma(E) 
\end{equation}
and
\begin{equation}
\labell{eq:Leibn22}
\{f, s\}(p)=0\quad\text{for all}\quad s\in \Gamma(E), 
\end{equation}
whenever $f$ vanishes at $p$ up to second order.
\end{Proposition}

\begin{Example}
\labell{rmk:loc-triv}
Let $E$ be the trivial line bundle $P\times \R^k$ and let us identify 
$\Gamma(E)$ with $C^\infty(P,\R^k)$. By setting the bracket
$\{~,~\}\colon \C^\infty(P)\times \Gamma(E)\to \Gamma(E)$ to be the 
component-wise Poisson bracket on $P$, we turn $E$ into a Poisson line bundle,
which we refer to as \emph{Poisson--trivial}. A Poisson vector bundle $E$ 
is said to 
be \emph{locally Poisson--trivial} if $E$ is locally isomorphic to a 
Poisson--trivial vector bundle. Poisson vector bundles are rarely
locally Poisson--trivial (see examples below).
\end{Example}

\begin{Remark}
A complex Poisson vector bundle $E$ is said to be \emph{Hermitian} or 
\emph{unitary} if $E$ is equipped with a Hermitian inner product
$\left<~,~\right>$ such that
$$
\{f,\left<s_1,s_2\right>\}=\left<\{f,s_1\},s_2\right>+
\left<s_1,\{\bar f,s_2\}\right>
$$
for all $f\in C^\infty(P)$ and all sections $s_1$ and $s_2$ of $E$.
In contrast with ordinary vector bundles, not every Poisson vector bundle
can be made Hermitian. The definition of Euclidean Poisson vector bundles is 
similar.
\end{Remark}

We emphasize that unless stated otherwise $\Vect_\pi$ means either one
out of four semi-rings: complex, real, Hermitian, or Euclidean Poisson
vector bundles. Of course, in every formula or equality vector bundles
or representations are assumed to be in the same class.

Let us illustrate these definitions by a series of examples.

\subsection{Basic examples}

\begin{Example}
\labell{exam:0}
Assume that the Poisson structure on $P$ is zero. Then every vector bundle
over $P$ with trivial bracket is a Poisson vector bundle.
\end{Example}

\begin{Example}
\labell{exam:flat}
Let $E$ be a bundle with a flat connection $\nabla$ over a Poisson
manifold $P$. Then $E$ is a Poisson bundle over $P$ with 
$\{f,s\}=\nabla_{X_f}s$. 
\end{Example}

\begin{Example}
\labell{exam:lie-alg}
Let $\g$ be a finite-dimensional Lie algebra and let $P=\g^*$ with the
standard Poisson structure. Let $V$ be a finite-dimensional representation
of $\g$. Set $E=\g^*\times V$ and $\{f,s\}(p)=df(p)\cdot s(p)+(L_{X_f} s)(p)$.
Here the first term denotes the action of $df(p)\in T_p^*\g^*=\g$ on
$s(p)\in V$ and the second term is just the Lie derivative of the vector valued
function $s$. It is easy to see that $E$ is a  
Poisson vector bundle over $\g^*$. This bundle is (locally) Poisson--trivial if and only
if $V$ is a trivial representation.
\end{Example}

The next two examples are more interesting from the Poisson geometry
perspective.

\begin{Example}
\labell{exam:Poisson-Lie}
Let $G$ be a Poisson Lie group and let $G^*$ be its dual, which we assume to be
simply connected. Then every representation of $G$ gives rise to a Poisson
vector bundle over $G^*$ similarly to Example \ref{exam:lie-alg}. We give
an explicit construction of this Poisson vector bundle in Example
\ref{exam:Poisson-Lie2}. The existence of such a Poisson
vector bundle is hardly surprising since $G^*$ is isomorphic to $\g^*$, when
$G$ is compact semisimple, \cite{Al,GW}.
\end{Example}

\begin{Example}
\labell{exam:bruhat}
Let $K$ be a compact Poisson Lie group with the standard Poisson structure
$\pi$ and let $T$ be the maximal torus in $K$ such that $\pi$ vanishes
along $T$. The Poisson structure $\pi$ descends to the quotient space 
$P=K/T$ making it into a Poisson manifold; see \cite{LW}. Let $V$ be a
representation of $T$. The vector bundle $E=K\times_T V$ over $P$ is
Poisson with the bracket defined as follows. Let us identify 
$\Gamma(E)$ with the space of $T$-equivariant functions $K\to V$ and
$C^\infty(P)$ with the space of $T$-invariant functions $K\to\C$. Then
$\{f,s\}$ is just the bracket of $f$ and $s$ thought of as functions on $K$,
i.e., $L_{X_f}s$, where $X_f$ is the Hamiltonian vector field of $f$ on $K$.
(It is not hard to see that $\{f,s\}$ so-defined is again $T$-equivariant
and hence a section of $E$.)
\end{Example}

The next two examples illustrate what happens when the requirement (L2)
is omitted. (See also Example \ref{exam:T*P} below.) 

\begin{Example}
\labell{exam:TP}
Let $E=TP$ and $\{f,s\}=[X_f,s]$, where $X_f$ is the Hamiltonian vector
field of $f$. This is \emph{not} a Poisson vector bundle unless $\pi=0$.
The bracket satisfies (L1), but not (L2): the condition
\eqref{eq:Leibn22} fails if $\pi\neq 0$.
\end{Example}

\begin{Example}
\labell{exam:pre-quant}
Let $P$ be symplectic with integral symplectic form $\omega$. Let 
$\LL\to P$ be a pre-quantum complex line bundle over $P$, i.e., a line
bundle with $c_1(\LL)=[\omega]$ and let $\nabla$ be a Hermitian connection 
on $\LL$ with curvature $-\omega$. Set $\{f,s\}=-\nabla_{X_f}s+\sqrt{-1}fs$. 
Then, $\LL$ is \emph{not} a Poisson vector bundle.
The bracket satisfies (L1) but not (L2): the condition \eqref{eq:Leibn21}
fails. The construction of pre-quantization generalizes to Poisson algebras 
associated with 
closed two-forms or Dirac structures (see \cite[Chapter 6]{GGK}) and 
to some Poisson manifolds, \cite{Hu,va:quant}. However, Poisson manifolds 
which are pre-quantizable in this sense are rather rare: the line bundle 
$\LL$ over $(P,\pi)$ exists if and only if $\pi$ is the image under 
$\pi^\#$ of an integral closed two-form on $P$, \cite{Hu,va:quant}. 
For example, the dual space of the Lie algebra of a compact semisimple Lie 
group is never pre-quantizable.
\end{Example}

\begin{Remark}
Example \ref{exam:pre-quant} shows that the class of vector bundles over
a Poisson manifold satisfying less restrictive conditions than (L1) and (L2)
may also be of interest. For instance, to include pre-quantization of Poisson
manifolds, one should keep (L1), but replace (L2) by the weaker condition
\eqref{eq:Leibn22}, which guarantees that $\{f,s\}(p)$ is determined by the 
1-jet of $f$ at $p$. It is also interesting to extend $\Vect_\pi(P)$
by adding to it representations up to homotopy, \cite{ELW}.
However, for our purposes the class of Poisson vector 
bundles is already too broad and we do not consider these larger semi-rings.
\end{Remark}

\begin{Remark}[Terminology]
\labell{rmk:term1}
We emphasize that neither the total space
of a Poisson vector bundle $E$ nor $E^*$ nor their fibers carry a natural 
Poisson structure. The condition that $E$ is a Poisson vector bundle 
can be expressed in terms a certain tensor field with values in $E$ and 
compatible with the Poisson structure on $P$ (cf., Section 
\ref{subsec:sec-ch}). However, here we prefer to
think of this condition as that $\Gamma(E)$ is a Poisson module over
$C^\infty(P)$.
\end{Remark}

\subsection{Poisson vector bundles as $\Omega^1(P)$-modules}
Recall that the space of one-forms $\df(P)$ on a Poisson manifold $(P,\pi)$
is a Lie algebra with the bracket
\begin{equation}
\labell{eq:bracket}
[\alpha,\beta]=-L_{\pi^\#\alpha}\beta+L_{\pi^\#\beta}\alpha
                 -d\pi(\alpha,\beta),
\end{equation}
where $\alpha$ and $\beta$ are in $\df(P)$. This bracket is a natural 
extension (via the graded Leibniz identity) of the bracket
on closed one-forms locally defined by $[df,dg]=d\{f,g\}$. 

For $\alpha\in\df(P)$ and $s\in\Gamma(E)$ we set 
\begin{equation}
\labell{eq:bracketE}
[\alpha,s](p)=\{f,s\}(p), 
\end{equation}
where $f\in C^\infty(P)$ is an arbitrary function function with 
$df(p)=\alpha_p$. From Proposition \ref{prop:Leibn2} it is easy to see that 
in this way we obtain a well-defined action of $\df(P)$ on $\Gamma(E)$.
The following result is an immediately consequence of the definition:

\begin{Proposition}
\labell{prop:bracket}
The bracket \eqref{eq:bracketE} turns $\sec$ into
a module over the Lie algebra $\df(P)$ such that $[\alpha,s]$ depends
only on $\alpha_p$ and the 1-jet of $s$ at $p$. Conversely, every vector
bundle equipped with a bracket $[~,~]\colon \df(P)\times \sec\to \sec$ 
satisfying these conditions is a Poisson vector bundle with
$\{f,s\}:=[df,s]$.
\end{Proposition}

The requirement that the bracket at $p$ depends only on $\alpha_p$ is equivalent
to the identity
\begin{equation}
\label{eq:bracketf}
[f\alpha,s]=f[\alpha,s]
\end{equation}
for all one-forms $\alpha$, sections $s$, and smooth functions $f$.

\begin{Example}[Continuation of Example \ref{exam:Poisson-Lie}] 
\labell{exam:Poisson-Lie2}
Let $G$ be a Poisson Lie group, $G^*$ its dual, and $V$ a representation of
the Lie algebra $\g$ of $G$. Let $E=G^*\times V$. Identifying 
$\alpha\in\g$ with a left-invariant one-form on $G^*$, we set
$[\alpha,s](p)=\alpha_p\cdot s(p)$ for $s\in\Gamma(E)$ and extend this
bracket to all one-forms on $G^*$ via \eqref{eq:bracketf}. As a result, $E$
becomes a Poisson vector bundle.
\end{Example}

Alternatively, we can think of the bracket $[\alpha,\cdot]$ as
a vector field on $E$:

\begin{Proposition}
\labell{prop:general}
For every $\alpha\in\df(P)$ there exists a unique vector field
$Z_\alpha$ on the total space of $E$ such that
\begin{enumerate}

\item The projection of $Z_\alpha$ to $P$ is equal to $\pi^\#\alpha$ and
hence the (local) flow $\varphi^t_\alpha$ of $Z_\alpha$ sends fibers of $E$ 
to fibers of $E$. Furthermore, $\varphi^t_\alpha$ is 
fiberwise linear and commutes with the action of the structural group 
of $E$ incorporated in the definition of $E$ as a Poisson vector 
bundle (e.g, $\U(n)$ if $E$ is Hermitian).

\item Let $Z_\alpha^*$ be the vector field induced by $Z_\alpha$ on $E^*$.
Then $L_{Z^*_\alpha}s=-[\alpha,s]$,
where $s$ on the left hand side is viewed as a fiberwise linear function on $E^*$.

\item $Z_{[\alpha,\beta]}=[Z_\alpha,Z_\beta]$ for all $\alpha$ and $\beta$
in $\df(P)$.

\item If $\alpha$ is closed, $\varphi^t_\alpha$ preserves the structure of 
a Poisson vector bundle.

\item $\varphi^t_\alpha$ is defined on the same interval
of time as the (local) flow of $\pi^\#\alpha$.
\end{enumerate}
\end{Proposition}

\begin{proof} 
The first and the second properties of $Z_\alpha$ comprise the
definition of this vector field. The third property is equivalent to
that $\Gamma(E)$ is a $\df(P)$-module; the fourth property is clear.
Finally, the last property follows from the first (see, e.g., the 
proof of \cite[Theorem 3.4]{Gi:mm}). We leave out the details of the
proof for the sake of brevity.
\end{proof}

\begin{Corollary}
\labell{cor:sympl}
Let $E$ be a Poisson vector bundle over a symplectic manifold $P$.
Then $E$ is locally trivial and, as in Example \ref{exam:flat},
$\{f,s\}=\nabla_{X_f}s$, where $\nabla$ is a flat connection on $E$. 
As a consequence, $\Vect_\pi(P)$ is in a one-to-one correspondence with 
$\Rep(\pi_1(P))$.
\end{Corollary}

\begin{Corollary}
\labell{cor:normal}
Denote by $N_p^*$ the co-normal Lie algebra to the symplectic leaf through
$p\in P$. (By definition, $N_p^*$ is the dual of the normal space to the 
leaf through $p$.) The bracket \eqref{eq:bracketE} restricts to a
well-defined representation of $N^*_p$ on $E_p$. For all $p$ in the
same leaf these representations are equivalent to one another.
\end{Corollary}

\subsection{Homogeneous Poisson vector bundles and homotopy equivalence}
The class of Poisson vector bundles is too broad for our purposes
as the next example indicates. (See also Example \ref{exam-ss}.)

\begin{Example}
\labell{exam:0gen}
Let $E$ be an ordinary vector bundle over a trivial Poisson manifold $P$ 
(i.e., $\pi=0$) and let $\sigma$ be a family, parameterized by $p\in P$, of 
representations of commutative Lie algebras $N_p^*=T_p^*P$ on $E_p$, i.e., 
$\sigma\in \Gamma(TP\otimes (E^*\otimes E))$. Then $\{f,s\}=\sigma(df) s$ 
or, equivalently, $[\alpha,s]=\sigma(\alpha)s$ 
makes $E$ into a Poisson vector bundle. 
\end{Example}

We propose two essentially different ways to reduce the class of Poisson
vector bundles: by imposing an additional constraint to eliminate undesirable
ones and by introducing a new equivalence relation which is weaker than
isomorphism.

\begin{Definition}
\labell{def:def}
Poisson vector bundles $E_0$ and $E_1$ over $P$ are \emph{deformation 
equivalent} or \emph{homotopic} if $E_0$ and $E_1$ are isomorphic to
vector bundles that can be connected by a 
family $E_t$ of Poisson vector bundles, smooth in $t\in [0,1]$. This is 
equivalent to requiring that $E_0$ and $E_1$ are isomorphic to the 
restrictions of some Poisson vector bundle $E$ over $P\times [0,1]$ to 
$P\times\{0\}$ and $P\times\{1\}$, respectively.
\end{Definition}

We denote the semi-ring of deformation equivalence classes of Poisson
vector bundles over $P$ by $\Vect_\pi^\defo(P)$. In the Hermitian or
Euclidean case, the vector bundles $E_t$ (or, equivalently, $E$) are required
to be in the same category as $E_0$ and $E_1$.

\begin{Definition}
\labell{def:hom}
A Poisson vector bundle $E$ is \emph{homogeneous} if every smooth
family $\psi_t\colon P\to P$ of Poisson diffeomorphisms with 
$\psi_0=\id$ admits a lift to a family $\tilde{\psi}_t$ of 
automorphisms of $E$ as a Poisson vector bundle with $\tilde{\psi}_0=\id$, 
i.e.,
\begin{equation}
\labell{eq:hom}
\tilde{\psi}_t^*\{f,s\}=\{\psi_t^* f, \tilde{\psi}_t^* s\}.
\end{equation}
(These lifts are not a part of the structure of a homogeneous
vector bundle.) This is equivalent to that for every Poisson 
vector field $X$ on $P$, there exists an $\R$-linear operator
$D_X\colon \Gamma(E)\to \Gamma(E)$ such that
\begin{equation}
\labell{eq:hom2}
D_X[\alpha,s]=[L_X\alpha,s]+[\alpha,D_Xs]
\quad\text{and}\quad D_X(fs)=(L_X f)s+f D_X s
\end{equation}
for all $s\in\Gamma(E)$, all $f\in C^\infty(P)$, and all $\alpha\in\df(P)$.
\end{Definition}

Note that when $\psi_t$ is given by a family of closed differential
one-forms, the existence of the lift $\tilde{\psi}_t$ is guaranteed by 
Proposition \ref{prop:general}. The essence of Definition \ref{def:hom}
is that a homogeneous
Poisson vector bundle over $P$ is as homogeneous as $P$ itself.
For example, when $E$ is homogeneous, the representations $E_p$ of $N_p^*$
are equivalent to one another for all points $p$ which can be obtained from
each other by families of Poisson diffeomorphisms. We denote the semi-ring
of isomorphism classes of homogeneous Poisson vector bundles over $P$
by $\Vect_\pi^\hom(P)$.

\begin{Example}
\labell{exam:sympl}
By Corollary \ref{cor:sympl}, every Poisson vector bundle over a symplectic 
manifold $P$ is homogeneous.
Two such vector bundles are homotopic if and only if the corresponding 
representations of $\pi_1(P)$ belong to the same path-connected component 
of $\Rep(\pi_1(P))$. 
Hence, $\Vect_\pi(P)=\Vect_\pi^\hom(P)=\Rep(\pi_1(P))$ and
$\Vect_\pi^\defo(P)=\pi_0(\Rep(\pi_1(P)))$.
\end{Example}

Let us now classify Poisson vector bundles over a Poisson manifold with 
trivial Poisson structure. 

\begin{Proposition}
\labell{prop:0}
Every Poisson vector bundle over a trivial Poisson manifold $P$
arises as described in Example \ref{exam:0gen} and all such vector
bundles are homotopic to one another: $\Vect_\pi^\defo(P)=\Vect(P)$. 
The only homogeneous Hermitian (or Euclidean) Poisson vector bundle 
over $P$ is the one with the trivial bracket $[~,~]=0$ as in 
Example \ref{exam:0}. 
\end{Proposition}

\begin{proof}
The first assertion is clear by Proposition \ref{prop:bracket}.
If $E$ is homogeneous, the representations $\sigma(p)$
are equivalent to each other for all $p\in P$. Furthermore, by applying 
diffeomorphisms fixing $p$, we conclude from \eqref{eq:hom} 
that the operators $\sigma_p(\alpha)\colon E_p\to E_p$, $\alpha\in T^*P$, are 
conjugate to each other for all $\alpha\neq 0$ and, hence, 
have the same eigenvalues. In particular,
the operators $\sigma_p(\alpha)$ and $\sigma_p(a\alpha)$
have the same eigenvalues for all constants $a\neq 0$. 
Since $\sigma_p(a\alpha)=a\sigma_p(\alpha)$, these eigenvalues must be zero.
This implies that $\sigma_p=0$, when $E$ is Hermitian or Euclidean. Hence, 
$[~,~]=0$.
\end{proof}

Proposition \ref{prop:0} does not extend to all homogeneous vector
bundles. In other words, for a Poisson manifold with zero Poisson structure
the semi-ring $\Vect^\hom_\pi(P)$ is in general bigger than $\Vect(P)$.

\begin{Example}
\labell{exam:nilp}
Let $P=\R$ and $E=\R\times \C^n$. Denote by $A$ the $n\times n$-Jordan block 
and set $[dx,s]=\sigma(dx)s=As$, where $x$ is the coordinate on $\R$. We
extend the bracket to all one-forms on $P$ by \eqref{eq:bracketf}. It is
easy to see that the resulting Poisson vector bundle is homogeneous.
\end{Example}

\begin{Proposition}
\labell{prop:product1}
Let $F$ be a simply connected symplectic manifold, $B$ a Poisson manifold,
and let $P=B\times F$ be equipped with the 
product Poisson structure. Then there is a one-to-one correspondence 
between Poisson vector bundles over $P$ and Poisson vector bundles over $B$.
The same holds for deformation equivalence classes. 
Thus $\Vect_\pi(P)\cong\Vect_\pi(B)$ and 
$\Vect_\pi^\defo(P)\cong\Vect_\pi^\defo(B)$. If the Poisson structure on $B$
is zero, $\Vect_\pi^\hom(P)\cong\Vect_\pi^\hom(B)$.
\end{Proposition}

\begin{proof} Proposition \ref{prop:product1} is a particular case of
Propositions \ref{prop:product-int} and \ref{prop:product-weak}.
Let us outline a direct proof, leaving
entirely aside the question of homogeneity.

Denote by $p\colon P=B\times F\to B$ the natural projection.
Fixing the direct product structure on $P$, we embed $C^\infty(F)$ into
$C^\infty(P)$ as a Poisson subalgebra. (Here we use the fact that the
symplectic structure on $F_y=p^{-1}(y)$ is independent of $y\in B$.)
Clearly, for a Poisson vector bundle $E$ over $B$, the pull-back $p^*E$ is a 
Poisson vector bundle over $P$, where the action of $\df(F)$ on $p^*E$
is given by the natural flat connection on $p^*E$ along the fibers of $p$.
It is easy to see that deformation equivalent vector bundles have 
deformation equivalent pull-backs.

Let us now construct the correspondence inverse to the pull-back.

Let $E$ be a Poisson vector bundle over $P$. The restriction
$E|_{F_y}$ to a fiber $F_y$ inherits the structure Poisson vector bundle. 
Namely, for $h\in C^\infty(F)$ and $s\in \Gamma(E|_{F_y})$, we set
$\{h,s\}=\{h,\bar s\}|_{F_y}$, where $\bar s$ is an arbitrary extension of
$s$ to a section of $E$. To verify that this bracket is well defined,
we need to show that if $s$ is a section of $E$ vanishing along $F_y$ and 
$h\in C^\infty(F)$, then $\{h,s\}$ vanishes along $F_y$. To this end note 
that locally $s=\sum f_j s_j$, where
$f_j$ are some functions on $P$ vanishing along $F_y$ and $s_j$ are
some sections. Hence, 
$$
\{h,s\}=\sum (\{h,f_j\}s_j+f_j\{h,s_j\}).
$$
Clearly, every term on the right hand side vanishes along $F_y$ and hence
so does the left hand side. 

By Corollary \ref{cor:sympl}, $E|_{F_y}$ is a trivial
vector bundle with the Poisson vector bundle structure arising from a
flat connection. In other words, we obtain a flat $p$-fiberwise connection
$\nabla$ on $E$. Let $E_B(y)$ be the space of flat sections of 
$E|_{F_y}$. The spaces $E_B(y)$ form a vector bundle $E_B$ over $B$. The
sections of $E_B$ can be identified with a subspace of $\Gamma(E)$ formed
by sections which $\nabla$-horizontal. It is easy to see that $\Gamma(E_B)$
is invariant under the bracket with elements of $p^* C^\infty(B)$ and thus
inherits the structure of a Poisson vector bundle via the embedding 
$p^*\colon C^\infty(B)\to C^\infty(P)$. Therefore,
$E=p^*E_B$. It is clear that the correspondence $E\mapsto E_B$ preserves
deformation equivalence. This completes the proof of the proposition.
\end{proof}

In Sections \ref{sec:mvb} and \ref{sec:m-w} we will prove some 
variations on the theme of Proposition 
\ref{prop:product1}. Here we mention only the most obvious one.

\begin{Proposition}
\labell{prop:product2}
Let $P$ be a Poisson manifold whose symplectic 
foliation is a fibration $P\to B$ with compact and simply connected 
fibers. Assume also that the cohomology class of the symplectic structure
on the fibers does not vary from fiber to fiber. Then there is a 
one-to-one correspondence between Poisson vector bundles over $P$ 
and over $B$. The same holds for homogeneous Poisson vector bundles or 
deformation equivalence classes. 
\end{Proposition}

\begin{proof}
The argument is similar to the proof of Proposition \ref{prop:product1} with
Moser's theorem used to introduce a local (in $B$) direct product structure
on $P$. 
\end{proof}

This result does not generalize to the case where the cohomology
class of the symplectic structure on $F$ varies from fiber to fiber;
see Example \ref{exam:non-product}.

\subsection{$K$-rings}
\labell{subsec:K-th}
\begin{Definition}
The \emph{homotopy Poisson $K$-ring}, $K^\defo_\pi(P)$, of a Poisson 
manifold $P$ is the Grothendieck rings associated with the semi-ring
$\Vect^\defo_\pi(P)$. The \emph{homogeneous Hermitian Poisson $K$-ring}, 
$KU^\hom_\pi(P)$, is the Grothendieck ring associated with the semi-ring
of homogeneous Hermitian Poisson vector bundles over $P$. The 
\emph{homogeneous orthogonal Poisson $K$-ring}, $KO^\hom_\pi(P)$, is 
defined similarly.
\end{Definition}

We emphasize that $K^\defo_\pi$ denotes one of the four deformation
$K$-rings: complex, real, Hermitian, or Euclidean. Likewise, $K$
is either complex (or equivalently, Hermitian) or real (equivalently,
orthogonal) ordinary $K^0$-functor. In the homogeneous case, we consider
only the Hermitian ring $KU^\hom_\pi$ and the Euclidean one $KO^\hom_\pi$.
In every formula, all $K$-rings are to be understood in the same
category. Note that none of the Poisson $K$-rings is
functorial with respect to Poisson maps.

From the results of the previous section we conclude that these $K$-rings
satisfy (K1) and (K2):

\begin{Corollary}
\labell{cor:k-first}
1. Let, as in Proposition \ref{prop:product1}, $F$ be a simply connected 
symplectic manifold, $B$ a Poisson manifold, 
and let $P=B\times F$ be equipped with the product Poisson structure. Then 
$K_\pi^\defo(P)=K_\pi^\defo(B)$ as rings. If $B$ carries the zero
Poisson structure, we also have $K_\pi^\defo(P)=K(B)$ 
and $KU_\pi^\hom(P)=KU(B)$.

2. Let $P$ be symplectic. Then $KU_\pi^\hom(P)=RU(\pi_1(P))$, where $RU$
stands for the ring of virtual unitary representations, and 
$K_\pi^\defo(P)$ is the Grothendieck ring generated by
 $\pi_0(\Rep(\pi_1(P)))$.
\end{Corollary}

Of course, similar results hold in the Euclidean case.

Example \ref{exam:nilp} readily implies that for the zero Poisson
structure the Grothendieck ring associated with $\Vect_\pi^\hom(P)$
can be strictly greater than $K(P)$, i.e., (K1) fails for this
ring. This is the reason that we do not ``nominate'' 
the Grothendieck ring of $\Vect_\pi^\hom(P)$ to
a Poisson $K$-ring. Nor do we
consider the Grothendieck ring associated with $\Vect_\pi(P)$, 
for this ring is enormous as Example \ref{exam:0gen} shows.

\section{Representations of Algebroids and Groupoids}
\labell{sec:algebroids}

The goal of this section is to illuminate the connection between 
Poisson vector bundles and representations of algebroids and groupoids.
We begin by recalling in Sections \ref{subsec:repr-alg} and 
\ref{subsec:repr-gr} the necessary definitions and known results important 
for what follows. We refer the reader to \cite{we:book,Ma} for a
general introduction and original references; additional references are
provided in the text as needed. In Section \ref{subsec:sec-ch} we outline
the construction of secondary characteristic classes of algebroid 
representations, due to Crainic, \cite{crainic}, in a way suitable for
our exposition and in Section \ref{subsec:ch-hom} 
we prove that for a homogeneous vector bundle these characteristic classes 
are $H^1_\pi(P)$-invariant.

\subsection{Representations of algebroids}
\labell{subsec:repr-alg}
Let $\A\to P$ be an algebroid over a manifold $P$
with anchor map $\kappa$ and bracket $[~,~]$.

\begin{Definition}
A \emph{representation} of $\A$ or an \emph{$\A$-equivariant vector bundle} 
$E$ over $P$ is a vector bundle $E$ equipped with
a bracket $[~,~]\colon \Gamma(\A)\times \Gamma(E)\to\Gamma(E)$ turning
$\Gamma(E)$ into a Lie algebra module over $\Gamma(A)$ and such that
\begin{equation}
\labell{eq:equiv}
\text{$[\alpha,s](p)$ depends only on $\alpha_p$}
\end{equation}
and
\begin{equation}
\labell{eq:lift}
[\alpha,fs]=(L_{\kappa(\alpha)}f)s+f[\alpha,s]
\end{equation}
for all sections $\alpha$ of $\A$, all sections $s$ of $E$, and all smooth
functions $f$ on $P$. 
\end{Definition}

Note that the condition \eqref{eq:equiv} is equivalent to
\eqref{eq:bracketf}. Equation \eqref{eq:lift} guarantees that $[\alpha,s](p)$ 
depends only on the 1-jet of $s$ at $p$ and expresses the fact
that the bracket $[~,~]$ is a lift of $\kappa$ to $E$ via infinitesimal
automorphisms of $E$.

\begin{Remark}[Terminology]
Calling $E$ a representation of $\A$ is consistent with regarding 
an algebroid as a generalization of a Lie algebra. The name an 
``$\A$-equivariant vector bundle'' comes from
a different perspective: we view an algebroid as a generalization of a Lie 
algebra action on a manifold; see Example \ref{exam:g-action}. The same 
applies to equivariant vector bundles over groupoids considered later on.

In the terminology of \cite{Fe1,Fe2}, representations of $\A$ are
vector bundles with flat connections over $\A$. 
\end{Remark}

\begin{Example}
\labell{exam:g-action}
Let $\A$ be the action algebroid $P\times \g$ for an infinitesimal action
of a Lie algebra $\g$ on $P$. Then a representation of $\A$ is
just a $\g$-equivariant vector bundle over $P$. (See, e.g., \cite{we:book}
for the definitions.)
\end{Example}

\begin{Example}
\labell{exam:fol}
Let $\A$ be the algebroid of vector fields tangent to a regular foliation
$\calF$ on $P$. Then a representation of $\A$ is a vector
bundle equipped with a partial flat connection along $\calF$. An example of
such a vector bundle is the normal bundle to $\calF$ with the Bott connection.
\end{Example}

\begin{Remark} 
\labell{rmk:prop-general}
Proposition \ref{prop:general}, with obvious changes of notations, holds 
for representations of algebroids.
\end{Remark}

\begin{Remark}[Holonomy, \cite{Fe2}]
\labell{rmk:holonomy}
Let $\gamma\colon [0,1]\to \A$ be a mapping with the property that
$\kappa(\gamma)=(p\circ\gamma)'$, where $p\colon\A\to P$ is the 
projection and the prime denotes the derivative with respect to time.
Following Fernandes \cite{Fe2}, we refer to such a path as an \emph{$\A$-path}.
There is a well-defined \emph{holonomy} map 
$h(\gamma)\colon E_{\gamma(0)}\to E_{\gamma(1)}$; see \cite{Fe2} (and 
also \cite{Fe1,GG} in the Poisson setting). This concept arises naturally
when we interpreter $E$ as a vector bundle with $\A$-connection. 
(Alternatively, pick a time-dependent section $\alpha_t$ of $\A$ such that 
$\alpha_t(p\gamma(t))=\gamma(t)$ and let, in the notation of
Proposition \ref{prop:general},  $\varphi_t$ be the time-dependent flow of
$Z_{\alpha_t}$. Then $h(\gamma)=\varphi^1|_{E_{\gamma(0)}}$; see \cite{GG}.)
Holonomy is defined in a setting 
more general than for representations of $\A$, but only for
representations it is in some sense homotopy invariant. We will use this
notion in the proof of Theorem \ref{thm:mvb}.
\end{Remark}

\begin{Example}
\labell{exam:T*P}
The original bracket $[~,~]$ on $\A$ does \emph{not} make $\A$ a 
representation of $\A$ on itself: the
bracket $[\alpha,s](p)$ of two sections $\alpha$ and $s$ of $\A$ depends on 
the values of $\alpha$ and $s$ at $p$. In other words, 
\eqref{eq:bracketf} and \eqref{eq:equiv}  fail. Likewise, when $P$ is a 
Poisson manifold, the bracket $\{f,s\}=[df,s]$, for $f\in C^\infty(P)$ and 
$s\in\Omega^1(P)$, satisfies (L1), but not (L2). 
\end{Example}

Recall that the algebroid $T^*P$ is canonically associated to 
every Poisson manifold $(P,\pi)$; see, e.g., \cite[Section 13.3]{we:book}. 
As a vector bundle, this algebroid is $T^*P$, the anchor map is $\pi^\#$, 
and the bracket is given by \eqref{eq:bracket}.
As an immediate consequence of Proposition \ref{prop:bracket}, we obtain

\begin{Corollary}
\labell{cor:poisson-equiv}
There is a one-to-one correspondence between 
Poisson vector bundles over $P$ and representations of $T^*P$.
\end{Corollary}

\begin{Example}[The canonical line bundle, \cite{ELW}]
\labell{exam:canonical-lb}
The canonical line bundle $\L=\wedge^nT^*P$, where $n=\dim P$, is
a Poisson vector bundle over $P$ with the action of $\alpha\in \df(P)$
on $s\in \Omega^n(P)$ given by 
$L_{\pi^{\#}(\alpha)}s+\left<\pi,d\alpha\right>s=\alpha\wedge (d\iota_{\pi}s)$.
\emph{The canonical line bundle is homogeneous.} (Set $D_X=L_X$ in 
\eqref{eq:hom2}.) Furthermore, $\L$ is trivial as a
Poisson line bundle if and only if $P$ is unimodular, i.e., $P$ admits
a volume form invariant under all Hamiltonian flows.

More generally, for an arbitrary algebroid $\A$, there is a canonical 
line bundle $Q_{\A}$ introduced in \cite{ELW}, which is also a representation 
of $\A$. These line bundles are important for the definition of the modular 
class.
\end{Example}

Clearly, the definition of deformation equivalence extends to
representations of algebroids and groupoids, considered in the next section.

\begin{Remark}
The notion of a homogeneous Poisson vector bundle does not
readily generalize to representations of algebroids. 
It appears to be essential that in Definition \ref{def:hom} $X$ is both
a one-cocycle on $T^*P$ and a vector field on $P$. Perhaps, the right
criterion for the ``correct'' generalization is to have Theorem 
\ref{thm:morita-weak} hold for the semi-ring of homogeneous vector bundles. 
\end{Remark}

\subsection{Representations of groupoids} 
\labell{subsec:repr-gr}

Let $\G$ be a smooth groupoid over a manifold $P$. 

\begin{Definition}
A \emph{representation} of $\G$ or a \emph{$\G$-equivariant vector bundle} 
over $P$ is a vector bundle $E$ 
equipped with linear maps $E_x\to E_y$ for all
$g\in\G$, where $x$ and $y$ are the source and, respectively, the target 
of $g$. These maps are assumed to depend smoothly on $g$.
\end{Definition}

\begin{Example}
\labell{exam:G-action}
Let $G$ be a Lie group acting on a manifold $P$ and let $\G=P\times G$
be the action groupoid. Similarly to Example \ref{exam:g-action},
representations of $\G$ are simply $G$-equivariant vector bundles 
over $P$.
\end{Example}

Let $\A$ be the algebroid of $\G$. Every representation of $\G$ 
is naturally a representation of $\A$. As in the case of actions of Lie 
groups, the converse is also true under a suitable simply connectedness 
hypothesis. The following result is a consequence of, e.g., 
\cite[Theorem 3.6]{MM}.

\begin{Proposition}
\labell{prop:integr}
Assume that the source--fibers of $\G$ are simply connected.
Then every representation of $\A$ integrates to
a representation of $\G$. As a consequence, there is a one-to-one 
correspondence between representations of $\A$ and $\G$.
\end{Proposition}

This proposition is not very hard to prove directly using, for
example, the notion of holonomy; see Remark \ref{rmk:holonomy}.

In what follows, we denote the semi-rings of isomorphism classes of 
representations of $\A$ or $\G$
by $\Vect_{\A}(P)$ and $\Vect_{\G}(P)$. Passing to infinitesimal
actions gives rise to a homomorphisms $\Vect_{\G}(P)\to\Vect_{\A}(P)$, which
is an isomorphism under the hypotheses of Proposition \ref{prop:integr}. The
resulting Grothendieck rings, for homotopy classes of representations, are 
denoted  by $K_{\A}^\defo(P)$ and $K_{\G}^\defo(P)$, respectively. We also
assume the other notation conventions of Section \ref{subsec:K-th}.

With a topological groupoid $\G$ over a manifold $P$, one can associate
the classifying space $B\G$. The construction of the space $B\G$ is
similar to Milnor's construction of the classifying space $BG$ for
a group $G$; see \cite{Bott, Segal}.

\begin{Example}
Let $\G=P\times G$ be the action groupoid for an action of a group $G$
on a manifold $P$. Then $B\G=(P\times EG)/G$, where $EG$ is a contractible
space with a free action of $G$ and $BG=EG/G$ is the classifying space
for $G$-principal bundles.
\end{Example}

An analogue of Borel's construction shows that every representation of 
$\G$ gives rise to a vector bundle over $B\G$. Thus we have
\begin{equation}
\labell{eq:BG}
\Vect_\G(P)\to \Vect(B\G),
\end{equation}
where $\Vect(B\G)$ is
the collection of isomorphism classes of ordinary vector bundles
over $B\G$. 

In general, the map \eqref{eq:BG} is not a bijection even when $G$ 
is compact; see \cite{Oliver}. However, as is the case with
actions of compact groups, for a Hausdorff proper groupoid 
the cohomology $H^*(B\G)$ is exactly the space 
where the characteristic classes of representations of $\G$
belong. In particular, the completion $H^{**}(B\G;\R)$ is the target
space for the $\G$-equivariant Chern character. With the exception of
some (important) examples, groupoids are rarely proper (and sometimes not
even Hausdorff), and Borel's construction we have just described
appears to be of little use when applied directly and naively to
non-proper groupoids. (See,
e.g., \cite{bott2} for more elaborate applications of Borel's 
construction in the non-proper case.) In the non-proper setting, the
corresponding $K$-ring is usually defined using non-commutative
geometry techniques; see \cite{connes-survey,connes-book} and also
\cite{He} for an excellent brief introduction. As we have
pointed out above, it is unlikely that $K$-rings considered in this
paper give meaningful results for seriously non-proper groupoids. (See, however,
Section \ref{sec:g}.) Yet it would still be interesting to see the exact 
relation between $K^\defo_\G$ and its ``non-commutative'' counterparts. 

In any case, for a Poisson manifold $P$, the algebroid $T^*P$ rarely
integrates to a groupoid and the constructions of this section apply
only in simplest examples. (See, e.g., Section \ref{sec:g}.)

\subsection{Secondary characteristic classes}
\labell{subsec:sec-ch}
The construction of secondary characteristic 
classes of representations of $\A$ is similar to that of 
characteristic classes of vector bundles with flat connections 
(see, e.g., \cite{KT}) or $\g$-structures, \cite{Fu}. 
In this section we recall the basic definition,
focusing on the case where $E$ is trivial, for only this
case is used in the subsequent sections. We refer the reader
to \cite{crainic} for a detailed treatment of the general case.
(See also \cite{Fe1,Fe2} for a relevant construction.)
We also restrict our attention to complex vector bundles. 
After suitable modifications, analogues of the results proved below hold 
for vector bundles with other  structural groups with $\gl(n)$ 
being replaced by the Lie algebra of the structural group.

Let $\A$ be an algebroid over $P$ and $E$ a complex representation of $\A$
which is trivial as an ordinary vector 
bundle. Fix a trivialization $E=P\times \C^n$. Then, for every
section $\alpha$ of $\A$, the vertical component of $Z_\alpha$ becomes a
$\gl(n)$-valued function on $P$ (see Proposition \ref{prop:general}
and Remark \ref{rmk:prop-general}). Thus, we obtain a $C^\infty(P)$-linear 
mapping $\Xi\colon\Gamma(\A)\to C^\infty(P,\gl(n))$, which we view as an
element in $\Gamma(\A^*)\times\gl(n)$. In other words,
\begin{equation}
\labell{eq:Xi}
[\alpha,s]=L_{\kappa(\alpha)}s+\Xi(\alpha)s,
\end{equation}
where we have identified $\Gamma(E)$ with $C^\infty(P,\C^n)$
using the trivialization. It is easy to see that $\Xi$
satisfies the Maurer--Cartan equation
\begin{equation}
\labell{eq:mc}
d_\A\Xi+\frac{1}{2}[\Xi,\Xi]_{\Lie}=0,
\end{equation}
where $d_\A$ is the differential in algebroid cohomology (see, e.g.,
\cite{we:book}) and as usual
$[\Xi,\Xi]_{\Lie}(\alpha,\beta)=2[\Xi(\alpha),\Xi(\beta)]$ in $\gl(n)$.
(Here and in what follows, we use the subscript ``$\Lie$'' to distinguish
this bracket from the Schouten bracket.)

\begin{Remark}
\labell{rmk:mc}
It is easy to see that every $\Xi$ satisfying the Maurer--Cartan equation
determines the structure of a Poisson vector bundle on a trivial vector 
bundle $E$.
\end{Remark}

From this remark we immediately obtain

\begin{Corollary}
\labell{cor:line-bundles}
The isomorphism classes of one-dimensional representations of $\A$ on the 
trivial line bundle are in a one-to-one correspondence with the first 
algebroid cohomology $H^1(\A;\R)$ for real or Hermitian representation and
with $H^1(\A;\C)$ for complex line bundles.
\end{Corollary}

For every cocycle $c$ of degree $k$ in the Chevalley--Eilenberg complex 
$C^*(\gl(n))$, we set $\Xi_c=c\circ \Xi\in \Gamma(\wedge^k\A^*)$. 
The proof of
the following proposition is standard; see, e.g., \cite{crainic,Fu,KT}.

\begin{Proposition}
\labell{prop:char-cl1}
The section $\Xi_c$ is $d_\A$-closed: $d_\A\Xi_c=0$.
The cohomology class $[\Xi_c]\in H^k(\A)$ in the algebroid cohomology of
$\A$ is determined on the cohomology class of $c$ in $H^*(\gl(n))$ 
This class depends of the trivialization of $E$. If $P$ is contractible,
$[\Xi_c]$ is independent of trivialization. Alternatively, when the 
trivializations used in the construction of $\Xi$ are unitary and 
$[c]\in H^*(\gl(n),\u(n))$, the class $[\Xi_c]$ is independent of 
trivialization.
\end{Proposition}
 
\begin{Remark}
\labell{rmk:independence}
The dependence of $\Xi_c$ on the trivialization is only up to the pull-backs
of the cohomology class of $c$ in $H^*(\GL(n))$ by all elements of
$[P,\GL(n)]$. This justifies the last assertion of the proposition since
the inclusion $\U(n)\hookrightarrow\GL(n)$ is a homotopy equivalence.
\end{Remark}

\begin{Example}
Assume that $P=\g^*$ and $E=\g^*\times \C^n$ is a Poisson vector bundle
corresponding to a representation of 
$\g$ in $\C^n$ as in Example \ref{exam:lie-alg}. Since $\g^*$ is
simply-connected, the class $[\Xi_c]$ is well defined. Let us
describe the characteristic class homomorphism $[c]\mapsto [\Xi_c]$ explicitly.
The representation
map $\g\to\gl(n)$ gives rise to a homomorphism $H^*(\gl(n))\to H^*(\g)$.
Then $[c]\mapsto [\Xi_c]$ factors as follows 
$$
H^*(\gl(n))\to H^*(\g)\to H^*(\g;C^\infty(\g^*))=H^*_\pi(\g^*).
$$
In the last identity we use the canonical identification
of the cohomology which holds already on the level of complexes; 
see, e.g., \cite{conn-smooth,lu}. The same holds for the action algebroid
on $\g^*$ associated with the coadjoint action.
\end{Example}

As usual, characteristic classes do not determine the representation.
For example, vanishing of all characteristic classes,
even for $H^*(\gl(n))$, does not guarantee that a representation is 
trivial; see, e.g., Example \ref{exam:nilp}. 
The secondary characteristic classes defined via elements
of $H^*(\gl(n),\u(n))$ carry even less information than 
those arising from $H^*(\gl(n))$. 
For instance, if the representation is unitary, all characteristic 
classes $[\Xi_c]$ for $[c]\in H^*(\gl(n),\u(n))$ vanish and, in fact,
these characteristic classes are some obstructions to making a representation
unitary.

\begin{Remark}[Rigidity, \cite{Fu}]
\labell{rmk:rigidity}
Most of the characteristic classes defined above are in fact invariants
of homotopy equivalence. For the sake of simplicity let us focus again
on trivial vector bundles. Denote by $\h$ the Lie algebra of the structural
group. Then $[\Xi_c]$ is an invariant of homotopy of $E$, provided that
$[c]\in \ker (\var\colon H^*(\h)\to H^{*-1}(\h;\h^*))$, \cite{Fu}. For
instance, in the Hermitian case, i.e., when $\h=\u(n)$, every class is
an invariant of homotopy. For $\h=\gl(n)$, the invariants of homotopy
are described as follows. Denote by $u_1,\ldots,u_n$ multiplicative generators
of $H^*(\gl(n))$ so that $\deg u_k=2k-1$. Then $[\Xi_c]$ is an invariant
of homotopy equivalence if and only if the factorization of $[c]$ in
these generators does not involve $u_1$; see \cite{Fu}.
\end{Remark}

\begin{Example}[\cite{crainic,ELW}]
\labell{exam:modular}
Let $\L$ be the canonical line bundle of a Poisson manifold $P$; see
Example \ref{exam:canonical-lb}. Then the characteristic class of
$\L$ corresponding to $c=\tr$ is the modular class of $P$. This class
is not rigid.
\end{Example}

\subsection{Secondary characteristic classes of homogeneous Poisson vector
bundles}
\labell{subsec:ch-hom} 
Let now $P$ be a Poisson manifold. A special feature of homogeneous Poisson 
vector bundles is that their secondary characteristic classes belong to a 
much smaller space than those of general Poisson vector bundles. 
Recall that the first Poisson cohomology $H^1_\pi(P)$ is a Lie algebra
(with the bracket induced by the Lie bracket of vector fields) and
$H^*_\pi(P)$ is a module over $H^1_\pi(P)$. 

\begin{Proposition}
\labell{prop:inv}
Assume that $E$ is homogeneous. Then for every cocycle $c$
the cohomology class $[\Xi_c]\in H^*_\pi(P)$ is $H^1_\pi(P)$-invariant.
\end{Proposition}

\begin{proof} In what follows we assume that the vector bundle $E$ is trivial.
The proof of the general case requires only superficial changes.
In the calculations below we use the sign convention of \cite{Fu}.

Fix a trivialization $E=P\times\C^n$ and thus identify sections of $E$
with functions $s\in C^\infty(P,\C^n)$. Let $X$ be a Poisson vector field on 
$P$. Since $E$ is homogeneous, there exists a function
$A\in C^\infty(P,\gl(n))$ such that $D_Xs=L_Xs+As$. 
Combining \eqref{eq:hom2} with \eqref{eq:Xi}, we obtain by a straightforward 
calculation that 
\begin{equation}
\labell{eq:der-Xi}
(L_X\Xi)(\alpha)=L_{\pi^{\#}\alpha}A+[\Xi(\alpha),A]_{\Lie}
\end{equation}
for any $\alpha\in\df(P)$.

Let now $c$ be a $k$-cocycle on $\gl(n)$. Consider the $(k-1)$-vector field
$b$ on $P$ defined by 
$$
b(\alpha_1,\ldots,\alpha_{k-1})=c(A, \Xi(\alpha_1), \ldots, \Xi(\alpha_{k-1})).
$$
We complete the proof by showing that $d_\pi b=-L_X\Xi_c$.
This equality is verified by the following (also straightforward) calculation.

First recall that 
\begin{eqnarray*}
(d_\pi b )(\alpha_1,\ldots, \alpha_k) 
&=& \sum_{i<j}(-1)^{i+j-1} 
b([\alpha_i,\alpha_j],
\alpha_1,\ldots,\widehat{\alpha_i},\ldots,\widehat{\alpha_j},\ldots,\alpha_k)\\
& & \quad+\sum_l(-1)^l
L_{v_l}b(\alpha_1,\ldots,\widehat{\alpha_l},\ldots,\alpha_k)\\
&=& \sum_{i<j}(-1)^{i+j-1} 
c(A, \Xi([\alpha_i,\alpha_j]),
\xi_1,\ldots,\widehat{\xi_i},\ldots,\widehat{\xi_j},\ldots, \xi_k)\\
& & \quad+\sum_l(-1)^l
L_{v_l}c(A,\xi_1,\ldots,\widehat{\xi_l},\ldots, \xi_k),
\end{eqnarray*}
where we set 
$$
v_m=\pi^{\#}\alpha_m\quad\text{and}\quad \xi_m=\Xi(\alpha_m)
\quad\text{for}\quad m=1,\ldots,k.
$$ 
The Maurer--Cartan equation \eqref{eq:mc} guarantees that 
$$
\Xi([\alpha_i,\alpha_j])=[\xi_i,\xi_j]_{\Lie}+L_{v_i}\xi_j-L_{v_j}\xi_i.
$$
Hence,
\begin{eqnarray}
\labell{eqn-db}
\qquad (d_\pi b )(\alpha_1,\ldots, \alpha_k) 
&=& \sum_{i<j}(-1)^{i+j-1} 
c(A, [\xi_i,\xi_j]_{\Lie},
\xi_1,\ldots,\widehat{\xi_i},\ldots,\widehat{\xi_j},\ldots, \xi_k)\\
\nonumber
& &\quad +\sum_{i<j}(-1)^{i+j-1} 
c(A, L_{v_i}\xi_j,
\xi_1,\ldots,\widehat{\xi_i},\ldots,\widehat{\xi_j},\ldots, \xi_k)\\
\nonumber
& &\quad -\sum_{i<j}(-1)^{i+j-1} 
c(A, L_{v_j}\xi_i,
\xi_1,\ldots,\widehat{\xi_i},\ldots,\widehat{\xi_j},\ldots, \xi_k)\\
\nonumber
& & \quad+\sum_l(-1)^l
L_{v_l}c(A,\xi_1,\ldots,\widehat{\xi_l},\ldots, \xi_k).
\end{eqnarray}
Obviously,
\begin{eqnarray*}
L_{v_l}c(A,\xi_1,\ldots,\widehat{\xi_l},\ldots, \xi_k)
&=& 
c(L_{v_l}A,\xi_1,\ldots,\widehat{\xi_l},\ldots, \xi_k)\\
& & \quad + \sum_{i<l}
(-1)^{i-1}c(A,L_{v_l}\xi_i,\xi_1,\ldots,\widehat{\xi_i},\ldots,
\widehat{\xi_l},\ldots, \xi_k)\\
& & \quad + \sum_{l<j}
(-1)^{j}c(A,L_{v_l}\xi_j,\xi_1,\ldots,\widehat{\xi_l},\ldots,
\widehat{\xi_j},\ldots, \xi_k).
\end{eqnarray*}
By plugging this into \eqref{eqn-db}, we eliminate the second and the 
third double sum. As a result of these cancelations, we obtain 
\begin{eqnarray*}
(d_\pi b )(\alpha_1,\ldots, \alpha_k) 
&=& \sum_{i<j}(-1)^{i+j-1} 
c(A, [\xi_i,\xi_j]_{\Lie},
\xi_1,\ldots,\widehat{\xi_i},\ldots,\widehat{\xi_j},\ldots, \xi_k)\\
& & \quad+\sum_l(-1)^l
c(L_{v_l}A,\xi_1,\ldots,\widehat{\xi_l},\ldots, \xi_k).
\end{eqnarray*}
By \eqref{eq:der-Xi}, we have $L_{v_l}A= (L_X\Xi)(\alpha_l)-[\xi_l,A]_{\Lie}$.
Therefore,
\begin{eqnarray*}
(d_\pi b )(\alpha_1,\ldots, \alpha_k) 
&=& \sum_{i<j}(-1)^{i+j-1} 
c(A, [\xi_i,\xi_j]_{\Lie},
\xi_1,\ldots,\widehat{\xi_i},\ldots,\widehat{\xi_j},\ldots, \xi_k)\\
& & \quad+\sum_l(-1)^{l-1}
c([A,\xi_l]_{\Lie},\xi_1,\ldots,\widehat{\xi_l},\ldots, \xi_k)\\
& & \quad+\sum_l(-1)^l
c((L_{X}\Xi)(\alpha_l),\xi_1,\ldots,\widehat{\xi_l},\ldots, \xi_k).
\end{eqnarray*}
The first two terms in this formula add up to zero because $c$ is
a cocycle on $\gl(n)$. Hence,
\begin{eqnarray*}
(d_\pi b )(\alpha_1,\ldots, \alpha_k) 
&= & \sum_l(-1)^l
c((L_{X}\Xi)(\alpha_l),\xi_1,\ldots,\widehat{\xi_l},\ldots, \xi_k)\\
&= & -(L_X \Xi_c)(\alpha_1,\ldots,\alpha_k).
\end{eqnarray*}
\end{proof}

\section{Poisson vector bundles over $\g^*$}
\labell{sec:g}
In this section we study Poisson vector bundles on the dual space
$\g^*$ of a Lie algebra $\g$ in connection with the requirement (K3).
In particular, we prove that every
such vector bundle is homotopic to the one associated with a representation
of $\g$ as in Example \ref{exam:lie-alg}. Furthermore, we prove some
results suggesting that for semisimple Lie algebras, every homogeneous
Hermitian Poisson vector bundle is isomorphic to the one associated with a 
representation.

\subsection{Classification up to homotopy}
Let $E\to \g^*$ be a Poisson vector bundle. Denote by $V$ the representation
of $\g=N_0^*$ on $E_0$ as in Corollary \ref{cor:normal}. Then as is shown
in Example \ref{exam:lie-alg}, the trivial bundle $\g^*\times V$ is a Poisson
vector bundle.

\begin{Theorem}
\labell{thm:homotopy}
The Poisson vector bundles $E$ and $\g^*\times V$ are homotopic.
\end{Theorem}

\begin{proof} The argument is by ``pushing the non-linearity of $E$ to
infinity''.
Pick a trivialization of $E$ so as to identify $E$ with $\g^*\times V$
as ordinary, not Poisson, vector bundles. With the trivialization fixed,
the structure of a Poisson vector bundle on $E$ is given by the
$\gl(n)$-valued vector field
$\Xi\colon \df(\g^*)\to C^\infty(\g^*,\gl(n))$  on $\g^*$ satisfying the 
Maurer--Cartan 
equation \eqref{eq:mc} as described in Section \ref{subsec:sec-ch}.

Using the structure of a linear space on $\g^*$, we write
$$
\Xi=\Xi_0+R,
$$
where $\Xi_0$ is the value of $\Xi$ at the origin. The constant $\gl(n)$-valued
vector field $\Xi_0$ on $\g^*$ also satisfies the Maurer--Cartan equation.
In fact, $\Xi_0$ determines the Poisson vector bundle $\g^*\times V$ associated
with the representation $V$. Note that the Maurer--Cartan equations for $\Xi$ 
and $\Xi_0$ imply that
\begin{equation}
\labell{eq:xi-R}
[\pi,R]+[\Xi_0,R]_{\Lie}+\frac{1}{2}[R,R]_{\Lie}=0,
\end{equation}
where $\pi$ is the Poisson structure on $\g^*$. 

To prove the theorem it suffices to connect
$\Xi=\Xi_1$ and $\Xi_0$ by a family of vector fields $\Xi_t$, $t\in [0,1]$,
satisfying the Maurer--Cartan equation (see Remark \ref{rmk:mc}).
Denote by $\varphi^t\colon \g^*\to \g^*$ the dilation by factor $t^{-1}$, 
i.e., $\varphi^t(x)=t^{-1}x$ for $x\in \g^*$ and $t\in (0,1]$. Let us show 
that $\Xi_t=t\varphi^t_*\Xi$ satisfies the Maurer--Cartan equation.

To this end we first observe that $\varphi^t_*\Xi_0=t^{-1}\Xi_0$ since $\Xi_0$
is a constant vector field.  
Hence $\Xi_t=\Xi_0+t\varphi^t_*R$.
Likewise, $\varphi^t_*\pi=t^{-1}\pi$, for $\pi$ is linear. Therefore,
\begin{eqnarray*}
[\pi,\Xi_t]+\frac{1}{2}[\Xi_t,\Xi_t]_{\Lie} 
& = & [\pi,\Xi_0+t\varphi^t_*R]+\frac{1}{2}[\Xi_0+t\varphi^t_*R,
\Xi_0+t\varphi^t_*R]_{\Lie}\\
&= & t[\pi,\varphi^t_*R]+t[\Xi_0,\varphi^t_*R]_{\Lie}
+\frac{t^{2}}{2}[\varphi^t_*R,\varphi^t_*R]_{\Lie}\\
& = & t^{2}[\varphi^t_*\pi,\varphi^t_*R]+
t^{2}[\varphi^t_*\Xi_0,\varphi^t_*R]_{\Lie}
+\frac{t^{2}}{2}[\varphi^t_*R,\varphi^t_*R]_{\Lie}\\
& = & t^{2}\varphi^{t}_*\Big(
[\pi,R]+[\Xi_0,R]_{\Lie}+\frac{1}{2}[R,R]_{\Lie}\Big)\\
& = & 0
\end{eqnarray*}
Here the first equality follows from that $\Xi_0$ satisfies the 
Maurer--Cartan equation and the last one is a consequence of \eqref{eq:xi-R}.

Finally, observe that $\Xi_t$ extends smoothly in $t$ to $t=0$ with value
$\Xi_0$ and hence $\Xi_t$, $t\in [0,1]$, is the required
homotopy between $\Xi=\Xi_1$ and $\Xi_0$.
\end{proof}

Clearly, Poisson vector bundles associated with representations in
the same path-connected component of $\Rep(\g)$ are homotopy equivalent. As a 
result, we conclude that the requirement (K3) holds for $K_\pi^\defo$:

\begin{Corollary}
For any Lie algebra, $\Vect^\defo_\pi(\g^*)=\pi_0(\Rep(\g))$ and, hence,
 $K_\pi^\defo(\g^*)$ is the Grothendieck ring  generated
by $\pi_0(\Rep(\g))$. In particular,  if $\g$ is semi-simple, 
$K_\pi^\defo(\g^*)=R(\g)$. 
\end{Corollary}

\subsection{Semisimple Lie algebras}
Throughout this section $G$ is a semisimple simply connected Lie group and 
$\g$ is the Lie algebra of $G$.

\begin{Theorem}
\labell{thm:compact}
Let $\g$ be semisimple of compact type and let $P\subset \g^*$ be a
$G$-invariant open set.\footnote{To circumvent the difficulties arising
in dealing with (equivariant) vector bundles over non-compact spaces, we assume
here that $P$ is (equivariantly) homotopy equivalent to a finite 
($G$-) CW-complex.
The same applies to other examples where (equivariant) vector bundles over
non-compact manifolds are considered.} Then there is a one-to-one 
correspondence between Poisson vector bundles over $P$ and $G$-equivariant 
vector bundles over $P$. Every Poisson vector bundle over $P$ is homogeneous.
Thus, $\Vect_\pi(P)=\Vect_\pi^\hom(P)=\Vect_\pi^\defo(P)=\Vect_G(P)$.
\end{Theorem}

\begin{proof}
The algebroid $T^*P$ is isomorphic to the action algebroid of the
coadjoint action on $P$. The latter integrates to the action algebroid
$P\times G$. By Proposition \ref{prop:integr}, there is a one-to-one
correspondence between Poisson and $G$-equivariant vector bundles over $P$.
As is shown in \cite{GW} (or in \cite{conn-smooth} for closed balls centered
at $0\in\g^*$),
every Poisson vector field on $P$ is Hamiltonian. Hence, a Poisson
vector bundle over $P$ is automatically homogeneous.
\end{proof}

\begin{Corollary}
\labell{cor:compact}
Under the hypotheses of Theorem \ref{thm:compact},
$KU_\pi^\hom(P)=K_\pi^\defo(P)=K_G(P)$, which is also equal to the
Grothendieck ring of $\Vect_\pi(P)$. (In the real case, one should
replace $KU_\pi^\hom$ by $KO_\pi^\hom$.)
\end{Corollary}

An interesting consequence of these results is that the Poisson $K$-rings
$KU_\pi^\hom$ and $K_\pi^\defo$ really depend on the Poisson structure, but
not only on the symplectic foliation.

\begin{Example}[Genuine dependence on Poisson structures]
\labell{exam:non-product}
Let $P=\su(2)^*\ssminus \{0\}$. In other words, $P=S^2\times (0,\infty)$
with the foliated symplectic structure depending linearly on $t\in (0,\infty)$.
By Corollary \ref{cor:compact}, $KU_\pi^\hom(P)=K_\pi^\defo(P)=
K_{\SU(2)}(S^2)$.  This ring is
$RU(\U(1))$, i.e., the ring of Laurent
polynomials in one variable of degree two. 
On the other hand, for the product Poisson structure on 
$S^2\times (0,\infty)$, which has the same symplectic foliation, 
$KU_\pi^\hom(P)=K_\pi^\defo(P)=K(\R)=\Z$, by
Corollary \ref{cor:k-first}.
\end{Example}

The question of classification of Poisson vector bundles over $\g^*$ for
a semisimple $\g$ of non-compact type is much more subtle. 
For instance, since every one-dimensional representation of $G$ is trivial,
an analogue of Theorem \ref{thm:compact} would imply that every 
Poisson line bundle over $\g^*$ is Poisson--trivial which is not true.

\begin{Example}
\labell{exam-ss} 
By Corollary \ref{cor:line-bundles}, Poisson line bundles 
over $\g^*$ are classified by $H^1_\pi(\g^*)=H^1(\g;C^\infty(\g^*))$. This 
cohomology is  infinite-dimensional, 
\cite{we:local}. Hence, the collection of isomorphism classes of such line 
bundles is infinite. Furthermore, combining this construction with that
of Example \ref{exam:nilp}, one should be able to construct homogeneous
Poisson vector bundles over $\g^*$ which do not correspond to representations.
However, as in Example \ref{exam:nilp}, these vector bundles are not
Hermitian.
\end{Example}

\begin{Conjecture}
\labell{conj:bundles} Let, as above, $\g$ be a semisimple Lie algebra of
non-compact type.

1. Every homogeneous Hermitian Poisson vector bundle over $\g^*$ is 
Poisson--trivial.

2. There is a one-to-one correspondence between germs of real analytic Poisson 
vector bundles over $\g^*$ and representations of $G$. 
\end{Conjecture}

This first part of this conjecture is essentially based on the lack of
examples. Note that since $G$ has no non-trivial finite--dimensional 
unitary representation, the conjecture implies that there is a  
one-to-one correspondence between homogeneous Hermitian Poisson vector 
bundles over $\g^*$ and finite--dimensional unitary
representations of $G$.
Here the only additional evidence is Proposition \ref{prop:inv} 
combined with another conjecture,
namely, that $H^*_\pi(\g^*)$ contains no non-zero $H^1_\pi(\g^*)$-invariant 
elements. 

The situation is more satisfactory in the real analytic setting. First recall
that the germ of a real analytic Poisson vector field on $\g^*$ is Hamiltonian
with a real analytic Hamiltonian, \cite{conn-analytic,conn-analytic-ad}.
By Corollary \ref{cor:line-bundles}, this implies that Conjecture
\ref{conj:bundles} holds for real analytic line bundles (real, Hermitian,
or complex). Furthermore, the following result confirms the conjecture
in the category of formal power series:

\begin{Theorem}
Every formal vector bundle over $\g^*$ is formally equivalent to the one
associated with a representation of $\g$.
\end{Theorem}

\begin{proof}
The argument is similar to the proof of Poincar\'e's theorem on the formal
normal form of a vector field and to Weinstein's proof of the formal 
linearization theorem for Poisson structures, \cite{we:local}. Fix a (formal) 
trivialization of $E$. We write $\Xi$ as a formal power
series on $\g^*$ and assume that this formal power series contains only
a zero order term and terms of degree greater than or equal to $r>0$:
$$
\Xi=\Xi_0+\Xi_r+\ldots .
$$
Here, as in the proof of Theorem \ref{thm:homotopy}, the zero
order term  $\Xi_0$ is associated with the representation of 
$\g=N^*_0$ on $E_0$ and
$\Xi_r$ is a $\gl(n)$-valued vector field on $\g^*$ whose coefficients
are homogeneous polynomials of degree $r$. As usual, the dots denote terms 
of degree higher than $r$.

The vector field $\Xi_r$ can be interpreted as a cochain on $\g$ with
coefficients in $S^r\g\otimes\gl(n)$. Denote by $d$ the differential
$C^1(\g;S^r\g\otimes\gl(n))\to C^2(\g;S^r\g\otimes\gl(n))$,
where $\gl(n)$ is turned into a $\g$-module via $\Xi_0$. 
The Maurer--Cartan equation \eqref{eq:mc} implies that
\begin{equation}
\labell{eq:cocycle}
d\Xi_r=d_\pi\Xi_r+[\Xi_r,\Xi_0]_{\Lie}=0,
\end{equation}
i.e., $\Xi_r$ is a cocycle.

Our goal is to find a formal change of trivialization $\Phi_r$ of $E$ 
which eliminates the $r$-th order term $\Xi_r$. 
We look for $\Phi_r$ in the form $\id + A_r$, where $A_r$ is a homogeneous
$\gl(n)$-valued polynomial on $\g^*$ of degree $r$. (In the Hermitian case,
$A_r$ is $\u(n)$-valued and $\Phi_r=\exp A_r$.) Under any formal change of 
trivialization $\Phi\colon \g^*\to\GL(n)$, the $\gl(n)$-valued vector field 
$\Xi$ transforms to $\Phi^*\Xi=\Phi\Xi\Phi^{-1}+(d_\pi\Phi)\Phi^{-1}$. Hence,
$$
\Phi_r^*\Xi=\Xi_0+(\Xi_r+[A_r,\Xi_0]_{\Lie}+d_\pi A_r)+\ldots .
$$
In other words, to eliminate $\Xi_r$, we need to find 
$A_r$ with $dA_r=-\Xi_r$. Since $\g$ is semisimple, 
$H^1(\g;S^r\g\otimes\gl(n))=0$ and, by \eqref{eq:cocycle}, $A_r$ does exist.

To eliminate all non-constant terms in the power series $\Xi$, we argue
inductively. Let $\Phi_1$ be a change of trivialization with $d A_1=-\Xi_1$
and $\Phi_2$ a change of trivialization eliminating the second order terms
in the resulting power series, etc. Then the composition 
$\Phi_r\circ\cdots\circ\Phi_1$ eliminates all terms up to degree $r+1$. 
Furthermore, the terms of any given degree in the expansion of
$\Phi_r\circ\cdots\circ\Phi_1$ stabilize as $r\to\infty$. Hence, 
$\Phi=\lim_{r\to\infty} \Phi_r\circ\cdots\circ\Phi_1$ is well defined as a 
formal power series. This is the required change of trivialization.
\end{proof}

\begin{Remark}
It appears plausible that by using the method of majorants as
in \cite{conn-analytic,conn-analytic-ad} one can control convergence
of the power series $\Phi$ and extend the theorem to
the real analytic case, thus proving the second part of Conjecture 
\ref{conj:bundles}.
\end{Remark}

\section{Morita Invariance of $K$-rings}
\labell{sec:mvb}
\subsection{Morita equivalence}
\labell{subsec:morita}
Recall that two Poisson manifolds $(P_1,\pi_1)$ and $(P_2,\pi_2)$ are
\emph{Morita equivalent} if there exists a symplectic manifold $W$
and submersions
\begin{equation}
\labell{eq:morita}
P_1\stackrel{\rho_1}{\leftarrow} W \stackrel{\rho_2}{\to} P_2
\end{equation}
with simply connected fibers such that $\rho_1$ is Poisson, $\rho_2$ 
is anti-Poisson, and the fibers of $\rho_1$ and $\rho_2$ are symplectic
orthogonal to each other, 
\cite{xu}.\footnote{The original definition Morita equivalence
given in \cite{xu} requires the projections $\rho_1$ and $\rho_2$
to be complete, \cite{we:local}. Following a suggestion of Alan Weinstein,
we omit this additional assumption, for
completeness is immaterial for the properties of Morita equivalence discussed 
below and for other results proved in this section, cf., 
\cite[Remark 4.7]{GG}.}
The essence of this definition is that a Morita equivalence
between Poisson manifolds gives rise to a bijection between their symplectic
leaves and corresponding leaves have equal first cohomology and anti-isomorphic
normal Poisson structures. 

Let us summarize the properties of Morita equivalence which are important for
what follows. (In (M1), except the last assertion, and (M2) it
suffices to assume that $\rho_1$- and $\rho_2$-fibers are connected.)

\begin{enumerate}
\item[(M1)] \cite{we:local}.
For a symplectic leaf $F\subset P_1$, the projection 
$\rho_2(\rho_1^{-1}(F))$ is a symplectic leaf in $P_2$. By symmetry, 
$F\mapsto \rho_2(\rho_1^{-1}(F))$ 
is a one-to-one correspondence between symplectic leaves of $P_1$ 
and $P_2$. The corresponding leaves have anti-isomorphic normal Poisson
structures and isomorphic first cohomology groups. 

\item[(M2)] \cite{we:local}.
The Poisson annihilator of $\rho^{*}_1 C^\infty(P_1)$ in $C^\infty(W)$
is $\rho^{*}_2 C^\infty(P_2)$, and vice versa. The manifolds $P_1$
and $P_2$ have equal spaces of Casimir functions, both isomorphic to
$\rho^{*}_1 C^\infty(P_1)\cap \rho^{*}_2 C^\infty(P_2)$.

\item[(M3)] \cite{gi-lu:morita}.
Morita equivalent Poisson manifolds $P_1$ and $P_2$ have isomorphic 
first Poisson cohomology spaces. More explicitly, there
is a natural isomorphism 
\begin{equation}
\labell{eq:morita-coh}
I_W\colon H^1_\pi(P_1) \stackrel{\cong}{\to}  H^1_\pi(P_2)
,\end{equation}
which is equal, up to a sign, to the one found in \cite{crainic}.
We will recall the definition of this isomorphism in Lemma 
\ref{lemma:gi-lu}. 
\end{enumerate}

\begin{Remark}
In spite of its name, Morita equivalence is \emph{not} an equivalence
relation. However, it becomes such when restricted the class of 
\emph{integrable} Poisson manifolds, i.e., those which admit global 
symplectic groupoids. (See \cite{xu} for more details.) For non-integrable
Poisson manifolds, the full meaning of Morita equivalence is not quite clear,
yet in what follows we make no integrability assumptions.
\end{Remark}

\subsection{Morita equivalence and $K$-rings}
\labell{subsec:mvb}
\begin{Theorem}
\labell{thm:mvb}
A Morita equivalence \eqref{eq:morita} induces an isomorphism between
semi-rings of Poisson vector bundles on $P_1$ and $P_2$, \cite{crainic},
homotopy equivalence classes of such vector bundles, and semi-rings of 
homogeneous vector bundles:
\begin{eqnarray}
\labell{eq:mvb1}
\Vect_\pi(P_1) &\cong& \Vect_\pi(P_2), \quad \text{\cite{crainic}},\\
\labell{eq:mvb2}
\Vect^\defo_\pi(P_1) &\cong& \Vect^\defo_\pi(P_2),\\
\labell{eq:mvb3}
\Vect^\hom_\pi(P_1) &\cong& \Vect^\hom_\pi(P_2).
\end{eqnarray}
\end{Theorem}

\begin{Corollary}
$K_\pi^\defo(P_1)=K_\pi^\defo(P_2)$ and $KU_\pi^\hom(P_1)=KU_\pi^\hom(P_2)$.
\end{Corollary}

The proof of \eqref{eq:mvb1} is given in \cite{crainic} in the integrable case
and a more general infinitesimal argument is also outlined there. Then 
\eqref{eq:mvb2} is an immediate consequence of the proof. Let us recall
some elements of the proof which are needed to verify
\eqref{eq:mvb3} and also used in Section \ref{sec:m-w}.

We begin by recalling the notion of an \emph{inverse-image} or 
\emph{pull-back} or \emph{induced algebroid}, \cite[p. 202-203]{HM}.
Let $\rho\colon Q\to P$ be a 
submersion and let $\A$ be an algebroid over $P$ with anchor $\kappa$. 
As a vector bundle, 
the inverse-image algebroid $\rho^{\star}\A$ is a subbundle in
$\rho^*\A\oplus TQ$ formed by pairs $(\alpha,X)$ with 
$\kappa(\alpha)=\rho_*X$
(We reserve the notation $\rho^{\star}\A$ for the inverse-image
algebroid to distinguish it from the pull-back vector bundle $\rho^*\A$.)
Note that $\ker\rho_*$ is a subalgebroid in $\rho^{\star}\A$ and 
$\rho^{\star}\A\cong\rho^*\A\oplus \ker\rho_*$ as vector bundles.
The bracket on $\rho^{\star}\A$ is given by
$$
[(f\alpha,X),(g\beta,Y)]=(fg[\alpha,\beta]+(L_Xf)\alpha-(L_Yg)\beta, [X,Y]),
$$
where $f$ and $g$ are smooth functions on $Q$, $\alpha$ and $\beta$ are 
pull-backs of sections of $\A$, and $X$ and $Y$ are vector fields on $Q$.
The anchor map on $\rho^{\star}\A$ is simply $(\alpha,X)\mapsto X$.

\begin{Example}
\labell{exam:pull-back}
The pull-back $\rho^{\star} TP$ is naturally isomorphic to $TQ$ with the 
isomorphism being the anchor map. This is easy to check in local
coordinates.
\end{Example}

For a representation $E$ of $\A$, its pull-back $\rho^*E$
is automatically a representation of $\rho^{\star}\A$ such that, in the 
notations of the above formula,
$$
[(\alpha,X), fs]=(L_Xf)s+f[\alpha,s],
$$
where $s$ is the pull-back of a section of $E$, \cite{HM}.

\begin{Lemma}[\cite{crainic}] 
\labell{lemma:cr}
Let the fibers of $\rho$ be $l$-connected
with $l\geq 1$. Then the pull-back $E\mapsto \rho^*E$ induces a bijection
\begin{equation}
\labell{eq:bijection}
\Vect_\A(P)\stackrel{\cong}{\to}\Vect_{\rho^{\star}\A}(Q). 
\end{equation}
Furthermore, the natural pull-back homomorphism 
$\rho^*\colon H^k(\A)\to H^k(\rho^{\star}\A)$ is an isomorphism for $k\leq l$.
\end{Lemma}

The first assertion readily follows from
the fact that the $\ker\rho_*$-action on $\rho^*E$ coincides with the one
arising from the natural flat connection on $\rho^*E$ along $\rho$-fibers.
The second assertion is proved in \cite{crainic}
by using an analogue of the Serre--Hochschild
spectral sequence (see, e.g., \cite{Fu}) for algebroids, with the subalgebroid
$\ker\rho_*$.

\begin{proof}[Proof of Theorem \ref{thm:mvb}] 
Let us set $\A_1=T^*P_1$
and $\A_2=T^*P_2$. By \eqref{eq:bijection}, to establish \eqref{eq:mvb1},
it suffices to show that
$\rho_1^{\star}\A_1$ and $\rho_2^{\star}\A_2$ are anti-isomorphic.
The required anti-isomorphism 
$J\colon \rho_1^{\star}\A_1\to\rho_2^{\star}\A_2$ is 
\begin{equation}
\labell{eq:id-alg}
J:(\alpha,X)\mapsto (\iota_X\omega-\alpha,X), 
\end{equation}
where $\omega$ is
the symplectic structure on $W$. Combined with \eqref{eq:bijection}, this
proves \eqref{eq:mvb1}. Moreover, it is clear that the bijection
\eqref{eq:bijection} induces a bijection on the level of homotopy equivalence 
classes thus giving \eqref{eq:mvb2}. Furthermore, \eqref{eq:mvb1} preserves 
the classes of Hermitian or Euclidean vector bundles.

To prove \eqref{eq:mvb3}, we will show that the bijection \eqref{eq:mvb1} sends
homogeneous vector bundles to homogeneous ones. We need the following fact, 
which is implicitly contained in \cite{gi-lu:morita} and explicitly in 
\cite{GG}:

\begin{Lemma}[\cite{GG,gi-lu:morita}]
\labell{lemma:gi-lu}
Let $\xi_1$ be Poisson vector fields on $P_1$. Then there exists a Hamiltonian
vector field $\xi$ on $W$ and a Poisson vector field $\xi_2$ on $P_2$ such
that 
$$
(\rho_1)_*\xi=\xi_1
\quad\hbox{\rm and}\quad
(\rho_2)_*\xi=-\xi_2
.
$$
Furthermore, $I_W([\xi_1])=[\xi_2]$.
\end{Lemma}

Let $E^{(1)}$ be a homogeneous Poisson vector bundle over $P_1$ and let 
$E^{(2)}$ be a Poisson vector bundle over $P_2$ which corresponds to
$E^{(1)}$ under \eqref{eq:mvb1}. Denote by $E$ the vector bundle
$\rho_1^*E^{(1)}=\rho_2^*E^{(2)}$ over $W$.

Let $\xi_1$ be a 
Poisson vector field on $P_1$. Denote by $\psi_1$ the (local) flow 
of $\xi_1$ in some time $t$, by $\tilde{\psi}_1$ its lift to $E^{(1)}$, and by
$\psi$ the (local) flow of $\xi$ in time $t$ on $W$. The flow
$\psi$ lifts to $E$. Indeed, for $x$ and $y$ in $W$ such that 
$\psi(x)=y$, we set

$$
\tilde{\psi}\colon  E_x\cong E_{\rho_1(x)}^{(1)}
\stackrel{\tilde{\psi}_1}{\to}
E_{\rho_1(y)}^{(1)}\cong E_y
.$$
(Here we use the fact that $\rho_1\psi=\psi_1$.) It is easy to see that 
$\tilde{\psi}$ is indeed a lift of $\psi$ via automorphisms of the
representation $E$ of $\rho_1^{\star}\A_1$, which covers, in the obvious 
sense, the lift $\tilde{\psi}_1$.

By Lemma \ref{lemma:gi-lu}, $\psi$ projects to the flow $\psi_2$ of
$-\xi_2$ on $P_2$. To complete the proof of \eqref{eq:mvb3} 
we only need to show that $\tilde{\psi}$ has a well-defined projection 
$\tilde{\psi}_2$ which is a lift of $\psi_2$ to $E^{(2)}$.
Let $x$ and $x'$ be two points in the same $\rho_2$-fiber. Then
for $y=\psi(x)$ and $y'=\psi(x')$ we also have $\rho_2(y)=\rho_2(y')$.
Our goal is to show that the following diagram is commutative:
$$
\begin{CD}
E_x @>\tilde{\psi}>> E_y\\
@Vh_xVV           @Vh_yVV\\
E_{x'} @>\tilde{\psi}>> E_{y'}
\end{CD}
$$
where the vertical arrows are the identifications 
$$
h_x\colon E_{x}\cong E_{\rho_2(x)}^{(2)}\cong E_{x'}
\quad\text{and}\quad
h_y\colon E_{y}\cong E_{\rho_2(y)}^{(2)}\cong E_{y'}.
$$

These identifications can also be interpreted in terms of holonomy. Connect
$x$ and $x'$ by a path $\gamma$ in the $\rho_2$-fiber. Then $h_x$ is just
the holonomy along $\gamma$ with respect to the natural flat 
$\rho_2$-fiberwise connection on $E$. On the other hand, $\gamma$ can also
be thought of as a $\ker(\rho_2)_*$-path and hence as a 
$\rho^{\star}_2\A_2$-path (see Remark \ref{rmk:holonomy}). Moreover, the
holonomy $h_x$ is equal to the holonomy $h(\gamma)$ from Remark 
\ref{rmk:holonomy}. Similarly, $h_y=h(\psi(\gamma))$.

Since $\tilde{\psi}$ is an automorphism of $E$ as of a 
$\rho_1^{\star}\A_1$-representation, $\tilde{\psi}$ commutes with holonomy
along $\rho_1^{\star}\A_1$-paths. The  anti-isomorphism $J$ defined by 
\eqref{eq:id-alg} interchanges (up to a sign) the structures of 
$\rho_1^{\star}\A_1$- and 
$\rho_2^{\star}\A_2$-representations on $E$. Therefore, $\tilde{\psi}$
also commutes with holonomy along $\rho_2^{\star}\A_2$-paths.
In particular, $h(\psi(\gamma))\tilde{\psi}=\tilde{\psi}h(\gamma)$ which
completes the proof.
\end{proof}

\begin{Remark}[Morita invariance of the modular class]
\labell{rmk:modular}
One can show that the canonical line bundle (see 
Example \ref{exam:canonical-lb}) of $P_1$ goes to the canonical line bundle
of $P_2$ under the correspondence \eqref{eq:mvb1}; see \cite{crainic}.
Passing to the first
characteristic class turns the correspondence \eqref{eq:mvb1} into 
the isomorphism $I_W$. This shows that $I_W$ sends the modular class of
$P_1$ to the modular class of $P_2$, \cite{crainic}:
\begin{equation}
\labell{eq:modular}
I_W(\modd (P_1))=\modd (P_2)
\end{equation}
in the self-explanatory notation of \cite{GG}.

Let us outline a direct proof \eqref{eq:modular}. Let $\mu_1$ and $\mu_2$
be volume forms on $P_1$ and $P_2$, respectively. Denote by $D$ the symplectic
$*$-operator on $W$, \cite{Br}. The forms $D\rho_1^*\mu_1$ and 
$\rho_2^*\mu_2$ have the same degree $\dim P_2$ and the same
null-space $\ker (\rho_2)_*$. Hence, 
$$
D\rho_1^*\mu_1=\pm e^H \rho_2^*\mu_2
$$
for some smooth function $H$ on $W$. One can easily verify
that $(\rho_1)_*X_H=v_1$ and $(\rho_2)_*X_H=-v_2$, where
$X_H$ is the Hamiltonian vector field of $H$ on $W$ and $v_1$ and $v_2$ are 
modular vector fields corresponding to $\mu_1$ and, respectively, $\mu_2$
(see \cite{we:modular}). 
By Lemma \ref{lemma:gi-lu}, $I_W([v_1])=[v_2]$ and \eqref{eq:modular}
follows.
\end{Remark}

\begin{Proposition}
\labell{prop:product-int}
Let $P$ be an integrable Poisson manifold and $F$ a simply connected
symplectic manifold. Equip $P\times F$ with the product Poisson 
structure. Then $\Vect_\pi^\hom(P)=\Vect_\pi^\hom(P\times F)$ and 
$KU_\pi^\hom(P)=KU_\pi^\hom(P\times F)$.
\end{Proposition}

\begin{proof} 
Since $P$ is integrable, $P$ is Morita equivalent to
$P\times F$, \cite{xu}.
\end{proof}

Of course, the assertion of Proposition \ref{prop:product-int} holds also
for $K_\pi^\defo$. Moreover, in this case the integrability assumption
can be omitted; see Proposition \ref{prop:product-weak} below.

\section{Weak Morita Category}
\labell{sec:m-w}
\subsection{Weak Morita equivalence}
\labell{subsec:general-alg}
To extend the results of Section \ref{subsec:mvb} to algebroids we
adopt the following definition,
which is a result of ``axiomatization'' of the proof of \eqref{eq:mvb1}. 

\begin{Definition}
Let $P_1$ and $P_2$ be smooth manifolds equipped with algebroids
$\A_1$ and $\A_2$, respectively. We say that $(P_1,\A_1)$ and $(P_2,\A_2)$
are \emph{weakly Morita equivalent} if there exists a manifold $Q$, 
submersions
$$
P_1\stackrel{\rho_1}{\leftarrow} Q \stackrel{\rho_2}{\to} P_2
$$
with simply connected fibers, and an isomorphism between the pull-back
algebroids $\rho_1^{\star}\A_1$ and $\rho_2^{\star}\A_2$.
\end{Definition}

It is clear that for Poisson manifolds weak Morita equivalence is indeed
weaker than genuine Morita equivalence. Weak Morita equivalence
has some properties similar to (M1)-(M3). For example, 
by Lemma \ref{lemma:cr}, 
weakly Morita equivalent algebroids $(P_1,\A_2)$ and $(P_2,\A_2)$ 
have isomorphic the zeroth and first algebroid cohomology spaces. The 
algebroids $\rho_1^{\star}\A_1$ and $\rho_2^{\star}\A_2$ have the same 
foliation. Hence, $F\mapsto \rho_2(\rho_1^{-1}(F))$, where $F$
runs through the leaves of $\A_1$, is a bijection between the leaves of
$\A_1$ and $\A_2$. Moreover, the corresponding leaves have isomorphic
first cohomology groups and even isomorphic fundamental groups as we will
soon show. Finally, if $x_1\in P_1$ and $x_2\in P_2$ are such that
$x_1=\rho_1(z)$ and $x_2=\rho_2(z)$ for some $z\in Q$, then $\A_1$ at
$x_1$ and $\A_2$ at $x_2$ have isomorphic stabilizers. 

The following result is a trivial consequence of Lemma \ref{lemma:cr}.

\begin{Theorem}
\labell{thm:morita-weak}
Assume that $(P_1,\A_1)$ and $(P_2,\A_2)$ are weakly Morita equivalent.
Then $\Vect_{\A_1}(P_1)\cong\Vect_{\A_2}(P_2)$ and 
$\Vect_{\A_1}^\defo(P_1)\cong\Vect_{\A_2}^\defo(P_2)$.
\end{Theorem}

\begin{Corollary}
\labell{cor:morita-weak}
$K_{\A_1}^\defo(P_1)\cong K_{\A_2}^\defo(P_2)$.
\end{Corollary}

\begin{Remark}[The category of representations]
In the algebraic setting, Morita equivalent rings have equivalent
categories of modules. The same is true for $C^*$-algebras,
\cite{Ri}. For a fixed algebroid $\A$ over $P$, let $\calV(\A)$ be the 
category of representations of $\A$. It is easy to see that for a submersion
$\rho\colon Q\to P$ with simply connected fibers, the pull-back 
$E\mapsto \rho^{*}E$ is an equivalence of categories 
$\calV(\A)\stackrel{\cong}{\to}\calV(\rho^{\star}A)$. Hence, \emph{weakly 
Morita equivalent algebroids have equivalent categories of representations}.
It would be interesting to check whether the converse is true or not.

A dual approach to ``modules'' over Poisson manifolds is developed by
Ping Xu, \cite{xu,xu2}, in the context of Morita equivalence of groupoids. 
In particular, Xu showed that integrable Poisson manifolds which
are Morita equivalent have equivalent ``categories'' of symplectic realizations
and ``symplectic left modules'', \cite{xu,xu2}.

\end{Remark}

\begin{Proposition}
\labell{prop:product-weak}
Let $P$ be a Poisson manifold and $F$ a simply connected
symplectic manifold. Equip $P\times F$ with the product Poisson 
structure. Then  $P$ and $P\times F$ are weakly Morita equivalent.
Thus $\Vect_\pi(P)\cong\Vect_{\pi}(P\times F)$ and
$\Vect_\pi^\defo(P)\cong\Vect_{\pi}^\defo(P\times F)$. As a consequence,
$K_\pi^\defo(P)=K_\pi^\defo(P\times F)$.
\end{Proposition}

\begin{proof} 
Let $\rho$ be the natural projection $P\times F\to P$. 
It suffices to show that $\rho^{\star}T^*P$ is isomorphic
to $T^*(P\times F)$. In the notation of the proof of Theorem \ref{thm:mvb},
the isomorphism is $(\alpha,X)\mapsto \alpha+\iota_{X_F}\omega$.
Here $X_F$ is the $F$-component of $X$, the form $\omega$ is the
symplectic form on $F$, and $\iota_{X_F}\omega$ is extended to a genuine
differential form on $P\times F$ using the direct product structure.
\end{proof}

Although weak Morita equivalence does not enable one to handle homogeneous
vector bundles, it is in some sense a ``better behaving'' notion
than Morita equivalence of Poisson manifolds:

\begin{Proposition}
\labell{prop:weak-me}
Weak Morita equivalence is an equivalence relation.
\end{Proposition}

We will prove this proposition in Section \ref{subsec:m-cat}.

Let us now consider some examples of weak Morita equivalence. The first
two of these are simply minor variations on the theme of the Morita equivalence
classification of symplectic manifolds by their fundamental groups;
see \cite[Proposition 2.1]{xu}.

\begin{Example}
\labell{exam:wm-TP}
Let $\A_1=TP_1$ and $\A_2=TP_2$. Assume that $P_1$ and $P_2$ are connected. 
The pairs $(P_1,TP_1)$ and $(P_2,TP_2)$ are weakly Morita equivalent 
if and only if $\pi_1(P_1)$ and $\pi_1(P_2)$ are isomorphic (cf.
\cite[Proposition 2.1]{xu}).

Indeed, as follows from Proposition \ref{prop:subm2}, weakly Morita 
equivalent manifolds always have isomorphic fundamental groups. To prove
the converse, set $Q=(\widetilde{P}_1\times \widetilde{P}_2)/\Gamma$. Here
$\widetilde{P}_1$ and $\widetilde{P}_2$ are universal coverings of $P_1$ and
$P_2$, respectively, and $\Gamma$ is the graph in 
$\pi_1(P_1)\times \pi_1(P_2)$ of an isomorphism $\pi_1(P_1)\cong \pi_1(P_2)$.
Clearly, the natural projections $\rho_1\colon Q\to P_1$ and 
$\rho_2\colon Q\to P_2$ have simply connected fibers. Furthermore,
by Example \ref{exam:pull-back}, we have isomorphisms
$\rho_1^{\star}TP_1\cong TQ\cong \rho_2^{\star}TP_2$.
\end{Example}

\begin{Example}
\labell{exam:wm-s}
Let $P_1$ and $P_2$ be connected symplectic manifolds with algebroids 
$T^*P_1$ and $T^*P_2$, respectively. The pairs $(P_1,T^*P_1)$ and 
$(P_2,T^*P_2)$ are weakly Morita equivalent if and only if $\pi_1(P_1)$ and 
$\pi_1(P_2)$ are isomorphic. By \cite[Proposition 2.1]{xu}, this is 
equivalent to that $P_1$ and $P_2$ are Morita equivalent. Indeed, for
weak Morita equivalence, algebroids matter only up to isomorphism.
On symplectic manifolds anchor maps are isomorphisms of algebroids:
$T^*P_1\cong TP_1$ and $T^*P_2\cong TP_2$. The assertion follows now from
Example \ref{exam:wm-TP}.
\end{Example}

\begin{Example}
\labell{exam:wm-0}
Let $P_1$ and $P_2$ be equipped with zero Poisson structures. (Hence, 
the algebroids are again $T^*P_1$ and $T^*P_2$, but both anchor maps are now
zero.) Then $P_1$ and $P_2$ are
weakly Morita equivalent if and only if these manifolds are diffeomorphic.
Indeed, in this case the submersions of a weak Morita equivalence 
$P_1\stackrel{\rho_1}{\leftarrow} Q \stackrel{\rho_2}{\to} P_2$
must have the same fibers.
\end{Example}

\begin{Example}
\labell{exam:wm-fol}
Let $P_1$ and $P_2$ carry foliations $\calF_1$ and $\calF_2$ and let
$\A_1=T\calF_1$ and $\A_2=T\calF_2$. Then $(P_1,\A_1)$ and $(P_2,\A_2)$
are weakly Morita equivalent if and only if there exist submersions
$P_1\stackrel{\rho_1}{\leftarrow} Q \stackrel{\rho_2}{\to} P_2$ with
simply connected fibers and such that the foliations $\rho_1^*\calF_1$
and $\rho_2^*\calF_2$ coincide.
\end{Example}

\subsection{Weak Morita category} 
\labell{subsec:m-cat}
Already the argument in \cite{gi-lu:morita}
suggested that Morita equivalence is an isomorphism in some category.
In this section, we introduce the corresponding category using weak Morita
equivalence which affords more flexibility than genuine Morita equivalence.

Fix a pair of integers $(k,l)$ with $0\leq k\leq l$. 
A \emph{weak Morita $(k,l)$-morphism} from an algebroid
$(P_1,\A_1)$ to an algebroid $(P_2,\A_2)$ is a pair of submersions
$P_1\stackrel{\rho_1}{\leftarrow} Q \stackrel{\rho_2}{\to} P_2$ such that
$\rho_1$-fibers are $l$-connected and $\rho_2$-fibers are $k$-connected
and an isomorphism of algebroids $\rho_1^{\star}\A_1\cong \rho_2^{\star}\A_2$.
(See \cite{HM} for the definition of  an isomorphism of algebroids.)

In what follows, we identify two morphisms between the same pairs
when there exists an diffeomorphism between their middle terms, $Q$ and $Q'$, 
and isomorphisms between the pull-back algebroids over $Q$ and $Q'$ such
that all the diagrams that arise are commutative. 

Let us define the composition of two weak Morita $(k,l)$-morphisms
$$
(P_1,\A_1)\stackrel{\rho_1}{\leftarrow} Q \stackrel{\sigma_1}{\to} (P,\A)
\quad\text{and}\quad
(P,\A)\stackrel{\sigma_2}{\leftarrow} S \stackrel{\rho_2}{\to} (P_2,\A_2).
$$
Let $W\subset Q\times S$ be the preimage of the diagonal under
$(\sigma_1,\sigma_2)\colon Q\times S\to P\times P$. It follows by 
transversality that $W$ is a smooth submanifold of $Q\times S$.
For $(x,y)\in W$
we set $\lambda_1(x,y)=\rho_1(x)$ and $\lambda_2(x,y)=\rho_2(y)$. 
It is easy to see that the composition is well-defined. (Recall
that a morphism is really an equivalence class.)

\begin{Theorem} 
\labell{prop:weak-mm}
$(P_1,\A_1)\stackrel{\lambda_1}{\leftarrow} W \stackrel{\lambda_2}{\to} 
(P_2,\A_2)$
is a weak Morita $(k,l)$-morphism with the isomorphism 
$\lambda_1^{\star}\A_1\cong \lambda_2^{\star}\A_2$ given
by \eqref{eq:mor} below.
\end{Theorem}

A routine verification, which we omit, now implies

\begin{Corollary}
The weak Morita $(k,l)$-category is indeed a category.
\end{Corollary}

Note that the identity map in this category is 
$P\stackrel{\id}{\leftarrow}P\stackrel{\id}{\rightarrow}P$.
We  call the pairs $(P_1,\A_1)$ and $(P_2,\A_2)$ \emph{weakly 
Morita $l$-equivalent} when there is a weak Morita $(l,l)$-morphism
between these pairs (but not necessarily an isomorphism). By symmetry,
weak Morita $l$-equivalence is an equivalence relation.

\begin{proof}[Proof of Proposition \ref{prop:weak-me}]
Weak Morita equivalence is just morphism in the weak Morita
$(1,1)$-category.
\end{proof}

\begin{proof}[Proof of Theorem \ref{prop:weak-mm}]
We focus on the properties of $\lambda_1$. For $\lambda_2$ the reasoning is 
identical.

To show that $\lambda_1$ is a submersion, we need to verify
that this map is onto and so is its derivative $(\lambda_1)_*$.
The surjectivity is a pure set-theoretic fact.
(Let $z\in P_1$. Take any $x\in\rho_1^{-1}(z)$ and any $y\in S$
with $\sigma_1(x)=\sigma_2(y)$. Then $\lambda_1(x,y)=z$.) This argument
uses only the assumption that all maps involved are onto. Hence, the
surjectivity also holds in the category of linear spaces and, in particular,
$(\lambda_1)_*$ is onto.

A similar argument shows that the natural projections
$Q\stackrel{p_1}{\leftarrow} W \stackrel{p_2}{\rightarrow} S$ are also 
submersions. 

Let us check that for $z\in P_1$, the fiber $\lambda_1^{-1}(z)$ is 
$l$-connected. The projection $\lambda_1^{-1}(z)\to\rho_1^{-1}(z)$ sending
$(x,y)$ to $x$ is a submersion whose fiber over $x$ is 
$\sigma^{-1}_2(\sigma_1(x))$. Hence, all fibers and the range of this 
submersion are $l$-connected. It suffices now to apply the following
proposition, which we prove in the next section.

\begin{Proposition}
\labell{lemma:subm}
Let $\rho\colon X\to Y$ be a submersion. Assume that $Y$ and all fibers of 
$\rho$ are $l$-connected. Then $X$ is also $l$-connected.
\end{Proposition}

It remains to define the isomorphism 
$\lambda_1^{\star}\A_1\cong\lambda_2^{\star}\A_2$.
First note the factorizations
$$
\lambda_1=\rho_1 p_1,
\quad
\lambda_2=\rho_2 p_2
\quad\text{and}\quad
\sigma_1 p_1=\sigma_2 p_2. 
$$
Using functoriality of the pull-back, \cite{HM}, we set this isomorphism to be
\begin{equation}
\labell{eq:mor}
\lambda^{\star}_1\A_1=p^\star_1\rho_1^{\star}\A_1\cong
p^\star_1\sigma_1^{\star}\A=p^\star_2\sigma_2^{\star}\A
\cong p^\star_2\rho_2^{\star}\A_2=\lambda_2^{\star}\A_2,
\end{equation}
where the isomorphisms are induced by the isomorphisms of the algebroids 
from the original weak Morita morphisms. This concludes the proof.
\end{proof}

By a straightforward inspection of the definition, we obtain

\begin{Proposition}
The algebroid cohomology up to degree $l$ and, for $l\geq 1$, 
$\Vect$, $\Vect^\defo$, and $K^\defo$ are contravariant functors on the weak 
Morita $(k,l)$-category. 
\end{Proposition}

In particular, weak Morita $(k,l)$-morphisms induce isomorphisms
of algebroid cohomology up to degree $k$ and also in $\Vect$, $\Vect^\defo$, 
and $K^\defo$ when $k\geq 1$.

\begin{Remark}
\labell{rmk:isom}
With the present equivalence of morphisms, the weak Morita
$(k,l)$-category has very few isomorphisms. A
$(k,l)$-morphism is an isomorphism if and only if both $\rho_1$
and $\rho_2$ are diffeomorphisms. (This can be seen by a simple
dimension count together with the description of $\lambda_1$ and
$\lambda_2$ in the proof of Theorem \ref{prop:weak-mm}.) For example,
consider the $(l,l)$-category. The morphisms 
$P_1\stackrel{\rho_1}{\leftarrow} Q \stackrel{\rho_2}{\to} P_2$ and
$P_2\stackrel{\rho_2}{\leftarrow} Q \stackrel{\rho_1}{\to} P_1$ are
inverse to each other if and only if $\rho_1$ and $\rho_2$ are
diffeomorphisms. This somewhat counterintuitive property suggests
that, perhaps, one should weaken the equivalence relation defining
morphisms.
\end{Remark}

Some examples of weak Morita $(k,l)$-morphisms can be readily
found by modifying Examples \ref{exam:wm-TP}-\ref{exam:wm-fol}.
The next three examples farther illustrate the naturality of various
constructions under weak Morita morphisms.

\begin{Example}[Weak Morita invariance of characteristic classes]
\labell{exam:wm-ch}
Consider the homomorphism of algebroid cohomology 
$H^*(\A_2)\to H^*(\A_1)$ induced by a weak Morita $(k,l)$-morphism from
$(P_1,\A_1)$ to $(P_2,\A_2)$. Then the Fernandes characteristic classes,
\cite{Fe2}, of $\A_2$ up to degree $l$ are mapped to the 
Fernandes characteristic classes of $\A_1$. This follows from
the fact that the pull-back $H^*(\A)\to H^*(\rho^{\star}\A)$, for a submersion
$\rho$, sends the characteristic classes of $\A$ to those of $\rho^{\star}\A$.

Additionally, as is easy to see, the canonical line bundle of $\A_1$ 
corresponds to the canonical line 
bundle of $\A_2$ (see Remarks \ref{exam:canonical-lb} and \ref{rmk:modular}).
\end{Example}

\begin{Example}[Weak Morita invariance of Poisson homology]
Let $P_1$ and $P_2$ be two weakly Morita $l$-equivalent
Poisson manifolds. Then \emph{the Poisson homology of $P_1$ and $P_2$ are 
isomorphic up to co-degree $l$, i.e., $H_{n_1-j}^\pi(P_1)=H_{n_2-j}^\pi(P_2)$ 
for $j\leq l$}, where $n_1=\dim P_1$ and $n_2=\dim P_2$. (See \cite{Br} for 
the definition of Poisson homology;
see also \cite{GG} for a much weaker version of this result.)
This follows from the duality between the Poisson homology and the cohomology
twisted by the canonical line bundle, \cite{ELW}, and from the fact that
the canonical line bundles of $P_1$ and $P_2$ are mapped to each other
by the equivalence. (The reader should be aware that the duality
is described in \cite{ELW} in different terms.) This example also generalizes
to algebroids.
\end{Example}

\begin{Example}
\labell{exam:homology}
A weak Morita $l$-equivalence induces isomorphisms of
homotopy and (co)homology groups up to degree $l$. This is a consequence of
Proposition \ref{prop:subm2}. 
\end{Example}

\begin{Remark}
The definition of weak Morita $(k,l)$-category is only one out of a few
inequivalent, but similar, constructions with the same properties. For 
instance, one may require the fibers to be simply connected and have 
vanishing real cohomology up to degree $k$ or $l$. All of the above 
results would still hold for this 
alternative definition. Another possible modification is to 
replace the isomorphism $\rho_1^{\star}\A_1\cong\rho_2^{\star}\A_2$ by
a morphisms of algebroids $\rho_1^{\star}\A_1\to\rho_2^{\star}\A_2$.
(See also Remark \ref{rmk:isom}.) 
\end{Remark}

This account of the weak Morita category is by no means intended to be
complete or detailed. Its only goal is to show that the concept is
natural and of some use and many questions are deliberately left out.
For instance: let $\calF_1$ on $P_1$ and $\calF_2$ on $P_2$ be 
Morita $l$-equivalent foliations. Does the isomorphism from
Example \ref{exam:homology} preserve their characteristic classes
(in the de Rham cohomology) up to degree $l$? 

\section{Proof of Proposition \ref{lemma:subm}} 
\labell{sec:proof}
Let us first set notations and terminology. In what follows,
$B$ and $K$ will be \emph{compact} manifolds possibly with boundary or
even with corners. (In fact, these will be either the sphere $S^n$,
the $n$-ball $D^n$, or the $n$-cube $I^n$.) For maps
$f\colon B\to Y$ and $F\colon B\times K\to X$, we say that $F$ 
\emph{covers} $f$ or that $F$ is a \emph{lift} of $f$ if $\rho(F(b,y))=f(b)$
for all $b\in B$ and $y\in K$, i.e., the diagram
$$
\begin{CD}
B\times K @>F>> X\\
@VVV @V{\rho}VV\\
B @>f>> Y
\end{CD}
$$
is commutative. The same applies to homotopies $f_t$ and $F_t$, where
$t$ is in some interval. 

Let us now establish some preliminary facts needed for the proof.

\begin{Lemma}[Local extensions of lifts]
\labell{lemma:local}
Let $f_t\colon B\to Y$, $t\in (-1,1)$, be a continuous homotopy and
$F_0\colon B\times K\to X$ a continuous lift of $f_0$. Then for a small 
$\eps>0$, there exists a continuous lift $F_t$ of $f_t$ for $t\in (-\eps,\eps)$.
\end{Lemma}

The lemma is clear: Since $F_0(B\times K)$ is compact, this set
has a neighborhood $U$ such that the restriction of $\rho$ to $U$ is
a trivial fiber bundle. 

\begin{Lemma}[Matching homotopies]
\labell{lemma:matching}
Let $f_t\colon B\to Y$, $t\in [-1,1]$, be a smooth homotopy.
Assume that there is
\begin{enumerate}
\item[(i)] a continuous lift $F_t^+\colon B\times K\to X$ of $f_t$ 
for $t\in [-1, 0]$,
\item[(ii)] a continuous lift $F_t^-\colon B\times K\to X$ of $f_t$ 
for $t\in [0,1]$
\item[] such that these two lifts coincide on $B\times \p K$ for $t=0$, i.e.,
$F_0^-(B\times \p K)=F_0^+(B\times \p K)$, and
\item[(iii)] a continuous homotopy $G_s$, $s\in [0,1]$,  relative 
$B\times \p K$ between $F_0^-$ and $F_0^+$ covering $f_0$. (We call $G$
a \emph{connecting homotopy}.)
\end{enumerate}
Then there exists a continuous lift $F_t$, $t\in [-1,1]$, of $f_t$ which 
is equal to $F_t^{\pm}$ for $t$ near $\pm 1$.
\end{Lemma}

We emphasize that this lemma is nearly obvious. In fact, by inserting $G$
between $F^-$ and $F^+$, we obtain a homotopy
which is almost a lift of $f$. This homotopy requires only a re-parameterization
and a minor alteration in the middle part to be turned into a genuine lift of
$f$. For the sake of completeness we give a detailed proof.

\begin{proof}
Assume first that $Y=B\times [-1, 1]$ and the homotopy $f$ is just the identity
map. Let
$$
Z=\bigcup_{s\in [0,1]} G_s(B\times K) \subset \rho^{-1}(B\times {0}).
$$
This is a compact set. It is easy to see that there exists a neighborhood
$U$ of $Z$ in $\rho^{-1}(B\times {0})$, a small $\eps>0$, and
an embedding $U\times (-\eps,\eps)\hookrightarrow X$ such that 
$\rho(x, t)=(\rho(x),t)$, where $x\in U$ and $t\in (-\eps,\eps)$. Here and in
what follows, we identify $U\times (-\eps,\eps)$ with a neighborhood of $Z$ 
in $X$. For small $|t|$, the homotopies $F_t^{\pm}$ take values in this
neighborhood.

It is also easy to see that $F_t^{\pm}$ can be changed for small $|t|$ so as
to make $F_t^{\pm}$ independent of $t$ in $U\times (-\eps,\eps)$. 
(In other words, the $U$-component of $F_t^{\pm}(y)$ is constant for
all $y\in B\times K$ and $|t|<\delta<\eps$.) From now on we assume that 
$F_t^{\pm}$ have this property. (The resulting maps can be chosen of
the same smoothness class as the original ones.)

Let $\phi\colon [-\delta,\delta]\to[0,1]$ be a smooth monotone
increasing function which is identically equal to $0$ for $t$ near
$-\delta$ and to $1$ for $t$ near $\delta$. The required lift is then
given by the formula
$$
F_t=\left\{
\begin{array}{l@{\quad\text{for}\quad}l}
F^-_t & t\in [-1,-\delta],\\
G_{\phi(t)} & t\in [-\delta,\delta],\\
F^+_t & t\in [\delta,1].
\end{array}
\right.
$$

Let us now prove the lemma for a general smooth homotopy
$f\colon B\times [-1,1]\to X$. Consider the pull-back submersion
$f^*\rho\colon \tilde{Y}\to B\times [-1,1]$ with the domain 
$$
\tilde{Y}=\big\{ (a,y)\in (B\times [-1,1])\times Y\mid \ f(a)=\rho(y)\big\}.
$$
By transversality, $\tilde{Y}$ is smooth and $f^*\rho\colon (a,y)\mapsto a$ 
is a submersion. A lift of $f$ pulls back to a continuous lift of $\id$.
(In fact, this pull-back is of the same smoothness class as the original lift, 
as one can check in local coordinates in the domain of $f$). 
Thus $F^+$, $F^-$, and $G$ pull back to maps satisfying
the hypotheses of the lemma and covering $\id$. For these-pull backs the
existence of the lift $F$ has already been proved. To obtain a lift of $f$,
it remains to compose this lift with the natural map $\tilde{Y}\to Y$.
\end{proof}

Now, with all of the preliminaries in place, we are in a position to
prove the key lemma.

\begin{Lemma}[Fillings of lifts]
\labell{lemma:key}
Let $B$ be the closed ball $D^n$ or the cube $I^n$ or the sphere $S^n$ 
and $K=S^m$. Let
$F\colon B\times K\to X$ be a continuous lift of a smooth map 
$f\colon B\to Y$. Then, for $n+m\leq l$, there exists a continuous lift 
$\Phi\colon B\times D^{m+1}\to X$ of $f$
(called a \emph{filling} of $F$) whose restriction to 
$B\times K= B\times \p D^{m+1}$ is $F$.
\end{Lemma}

We will need this lemma only for $m=0$, but it is more convenient to
prove it for all $m$ with $n+m\leq l$.

\begin{proof}[Proof of Lemma \ref{lemma:key} for $B=D^n$ and $I^n$]
We will prove the lemma for $B=I^n$. (This
implies the lemma for $B=D^n$, for $f$ can always be extended from
$D^n$ to some $n$-cube containing $D^n$ and $F$ can also be extended to
$I^n\times K$.) The proof is by induction in $n$.

For $n=0$, i.e., $I^0=\{\pt\}$, the assertion follows from the assumptions that 
the fibers of $\rho$ are $l$-connected and $m\leq l$. Assume that the
assertion is proved for $\dim B <n$. 

Set $I^n=I^{n-1}\times [0,1]$ and $C_t=I^{n-1}\times \{t\}$ for $t\in [0,1]$. 
Let also $f_t=f|_{C_t}$ and $F_t=F|_{C_t\times K}$.
By the induction hypothesis, $F_0$ has a filling $\ol{\Phi}_0$ and hence,
by Lemma \ref{lemma:local}, this filling extends over 
$I^{n-1}\times [0,\eps)$ for some
small $\eps>0$. Consider the maximal open in $[0,1]$ interval 
$[0,r)$ (or $[0,1]$) over which a filling $\ol{\Phi}_t$ of $F_t$ exists. 
If this interval
is $[0,1]$ the proof is finished. Hence, let us assume that the interval is
$[0, r)$ with $0<r\leq 1$ and arrive to a contradiction.
We will present the proof for $r<1$. The argument for $r=1$ requires only
minor modifications.

By the induction hypothesis, $F_r\colon C_r\times K\to X$ admits a filling,
which extends to a filling $\widetilde{\Phi}_t$ over some small interval 
$[r-\eps,r+\eps]$, by Lemma \ref{lemma:local}. Now we have two fillings,
$\ol{\Phi}_{r-\eps}$ and $\widetilde{\Phi}_{r-\eps}$, of $F_{r-\eps}$.
These are maps $C_{r-\eps}\times D^{m+1}\to X$ which coincide on 
$C_{r-\eps}\times S^{m}$ and hence give rise to a map 
$C_{r-\eps}\times S^{m+1} \to X$ covering $f_{r-\eps}$.\footnote{This map may
fail to be smooth at the points of $C_{r-\eps}\times S^m$ even if all 
the maps considered prior to this moment are.
Although a smoothing argument can be applied here, we prefer
to work with continuous lifts. In fact, this is the only
reason why we ever need to drop the smoothness assumption;
the proofs of Lemmas \ref{lemma:local} and \ref{lemma:matching} go through
in the smooth category.}

 The induction 
hypothesis applies to this map, for $C_{r-\eps}$ is a cube of
dimension $n-1$ and $(n-1)+ (m+1)=n+m\leq l$. Thus it has a filling,
which can be thought of as a homotopy $G$ between $\ol{\Phi}_{r-\eps}$ and 
$\widetilde{\Phi}_{r-\eps}$ which covers $f_{r-\eps}$.

Lemma \ref{lemma:matching} applies to 
the homotopy $\ol{\Phi}_t$ for $t\in [0,r-\eps]$, the homotopy 
$\widetilde{\Phi}_t$ for $t\in [r-\eps,r+\eps]$, and the connecting 
homotopy $G$. As a result, we obtain a filling of $F_t$ for 
$t\in [0, r+\eps]$. This contradicts the choice of $r$.
\end{proof}

\begin{proof}[Proof of Lemma \ref{lemma:key} for $B=S^n$]
As in the proof for $B=D^n$, we argue inductively. For $n=0$ and $m\leq l$,
the image $F(S^0\times S^m)$ is comprised of two mapped $m$-spheres contained
in the fibers of $\rho$. Since these fibers are $l$-connected, the filling
does exist. 

Assuming that the lemma is proved for $\dim B <n$, we establish it
for $S^n$. Let us think of $S^n$ as the union of two hemispheres:
$S^n=D^n_-\cup D^n_+$. Lemma \ref{lemma:key} has already been proved
for the discs and, hence, it applies to $F_{\pm}=F|_{D^n_{\pm}\times S^m}$.
As a result we obtain two fillings $\Phi_{\pm}$ over the hemispheres.
Along the ``equator'' $S^{n-1}\times S^m$ of $S^n\times S^m$, 
$\Phi_{\pm}$ are fillings of the same map $F|_{S^{n-1}\times S^m}$. Hence,
$\Phi_{\pm}$ give rise to a continuous map $S^{n-1}\times S^{m+1}\to X$ 
covering $f|_{S^{n-1}}$. By the induction hypothesis, this map has a filling 
and this filling can be viewed as a connecting homotopy $G$ between 
$\Phi_{-}|_{S^{n-1}\times S^m}$ and $\Phi_{+}|_{S^{n-1}\times S^m}$.
Applying Lemma \ref{lemma:matching} to
$\Phi_-$, $\Phi_+$, and $G$, we obtain the required filling for $F$.
\end{proof}

\begin{Lemma}[Lifting discs]
\labell{lemma:discs}
For every smooth map $f\colon D^n\to Y$, where $n\leq l+1$, there
exists a continuous lift $F\colon D^n\to X$.
\end{Lemma}

\begin{Remark}
Once a continuous lift is found, it is not hard show that a smooth lift
also exists. 
\end{Remark}

\begin{proof}[Proof of Lemma \ref{lemma:discs}]
We only outline the argument, for it is very similar to the proof 
of Lemma \ref{lemma:key}. As in that proof we replace $D^n$ by $I^n$.

Arguing inductively, we assume that the lemma is proved for the cubes of
dimension less than $n$. As above, we set $C_t=I^{n-1}\times \{t\}$
and $f_t=f|_{C_t}$, where
$I^n=I^{n-1}\times [0,1]$, and consider the maximal open interval in $[0,1]$ over
which the lift exists and which contains $0$. Denote this lift by $\ol{F}_t$.
Let us assume that the
interval has the form $[0,r)$ with $0<r\leq 1$ and arrive to contradiction.
As in the proof of Lemma \ref{lemma:key}, we focus on the case $r<1$. (The
case $r=1$ is handled similarly.)

By the induction assumption, we can lift $f_r$ to $X$ and then, using
Lemma \ref{lemma:local}, extend it to a homotopy 
$\widetilde{F}_t$, $t\in [r-\eps,r+\eps]$, for a small $\eps>0$. 
Thus we have two lifts of $f_{r-\eps}$, namely, 
$\ol{F}_{r-\eps}$ and $\widetilde{F}_{r-\eps}$. 
Together these lifts give rise to
a map $I^{n-1}\times S^0$ covering $f_{r-\eps}$. Applying Lemma
\ref{lemma:key}, we find a filling for this map. This filling 
is nothing else but a homotopy between $\ol{F}_{r-\eps}$ and 
$\widetilde{F}_{r-\eps}$ covering $f_{r-\eps}$. By Lemma \ref{lemma:matching},
the lift of $f_t$ extends to the interval $[0, r+\eps]$, which contradicts
the choice of $r$.
\end{proof}

\begin{proof}[Proof of Proposition \ref{lemma:subm}]
Let $F_0\colon S^n\to X$, where $n\leq l$, be a smooth map and let 
$f\colon D^{n+1}\to Y$ be a smooth homotopy of $\rho F_0=f|_{S^{n-1}}$ to
a point. We will find a continuous homotopy $F\colon D^{n+1}\to X$ of $F_0$
to a point. By Lemma \ref{lemma:discs}, there exists a lift
$\ol{F}\colon D^{n+1}\to X$ of $f$. Furthermore, by Lemma \ref{lemma:local}, 
$F_0$ can be extended to a little closed collar $U$ of $S^{n}$ in $D^{n+1}$.
Let $\widetilde{F}$ be this extension. Denote by $C$ the inner boundary
of $U$ in $D^{n+1}$. Thus $C$ is the sphere $S^n$ over which we have two lifts
$\ol{F}|_C$ and $\widetilde{F}|_C$ of $f|_C$. These lifts comprise a
map $C\times S^0\to X$ covering $f|_C$, which, by Lemma \ref{lemma:key}, 
have a filling. This filling is a homotopy $G$ between $\ol{F}|_C$ and 
$\widetilde{F}|_C$ covering $f|_C$. Applying Lemma \ref{lemma:matching},
we obtain the required lift $F$.
\end{proof}

\begin{Remark}
Proposition \ref{lemma:subm} does not hold if the assumption that $\rho$ is
a submersion is omitted and, as stated, Lemma \ref{lemma:discs} fails without
the assumption that $\rho$-fibers are $l$-connected.
\end{Remark}

A similar argument proves the following

\begin{Proposition}
\labell{prop:subm2}
Let $\rho\colon X\to Y$ be a submersion whose fibers are $l$-connected.
Then $\rho$ induces an isomorphism of homotopy groups
$\pi_n(X)\stackrel{\cong}{\to} \pi_n(Y)$ for $n\leq l$. 
\end{Proposition}

\end{document}